

\font\tenmsy=msym10
\font\sevenmsy=msym7
\font\fivemsy=msym5
\newfam\msyfam
\textfont\msyfam=\tenmsy  \scriptfont\msyfam=\sevenmsy
  \scriptscriptfont\msyfam=\fivemsy
\def\hexnumber@#1{\ifcase#1 0\or1\or2\or3\or4\or5\or6\or7\or8\or9\or
	A\or B\or C\or D\or E\or F\fi }
\edef\msy@{\hexnumber@\msyfam}
\mathchardef\ltimes="2\msy@6E

\centerline {{\bf ELLIPTIC ALGEBRAS AND}}
\vskip 1mm
\centerline {\bf EQUIVARIANT ELLIPTIC COHOMOLOGY I.}
\vskip 3mm
\centerline {({\it technical report})}

\vskip .7cm

\centerline {\bf V. Ginzburg, M. Kapranov, E. Vasserot}

\vskip 1cm

\centerline{Table of contents :}

\vskip .3cm

{$\qquad${\bf 1.} Equivariant elliptic cohomology : axiomatics.}

{$\qquad${\bf 2.} Thom sheaves and Gysin maps.}

{$\qquad${\bf 3.} Geometric construction of current algebras.}

{$\qquad${\bf 4.}  Current algebra on an elliptic curve and elliptic
cohomology.}

{$\qquad${\bf 5.}
Classical elliptic algebras and elliptic cohomology.}

\vskip 1cm

The motivation  for  writing this paper was the parallelism in the
classification of three different kinds of mathematical objects:

\vskip .1cm

\item{\bf (i)} Classical $r$-matrices.

\vskip .1cm

\item{\bf (ii)} Generalized cohomology theories that have Chern classes
for complex vector bundles.

\vskip .1cm

\item{\bf (iii)} 1-dimensional formal groups.

\vskip .2cm

\noindent  Recall first, that
there are some distinguished  1-dimensional formal groups - those corresponding
to actual algebraic groups - that is,
the additive group ${\bf G}_a$,
the multiplicative group ${\bf G}_m$, and the elliptic curves.
Further, classical $r$-matrices, that is
solutions to the classical Yang-Baxter equation
$$[r^{13}(u+v), r^{23}(v)-r^{12}(u)] = [r^{23}(v), r^{12}(u)]. \leqno (0.1)$$
were classified by Belavin and Drinfeld
[BD]. They showed
that , under appropriate non-degeneracy
conditions, all algebraic
solutions to (0.1) are
given by rational,
trigonometric and elliptic functions respectively. The ``quantum algebras"
associated to those solutions are known, respectively, as Yangians,
quantized enveloping algebras, and quantum elliptic algebras [Dr1]
[RS]. Thus, one gets an `a posteriori' correspondence
(i)$\leftrightarrow$(iii).

\vskip .2cm

The relation (ii)$\leftrightarrow$(iii),
on the other hand, is well known in
algebraic topology,
due to
Quillen [Q]. The generalized cohomology theories associated to formal
groups arising from
${\bf G}_a$, ${\bf G}_m$
and elliptic curves are, respectively, the ordinary cohomology,
the complex $K$-theory and the so-called elliptic cohomology  [Lan][S2].

\vskip .1cm
This suggests looking for a direct construction relating
one of the three generalized cohomology theories above with the corresponding
quantum algebras.
For the case of $K$-theory
such a relation was found in [GV1] following an earlier work [BLM].
It provides a description of a quantized universal enveloping algebra
in terms of equivariant
K-theory of flag varieties. It was also observed in
[GV1] that
replacing  equivariant K-theory by equivariant
homology yields a construction of
Yangians.

The aim of the present paper is to treat the elliptic case.
The first obstacle here is the absence of a `good'
equivariant elliptic
cohomology theory. The non-equivariant elliptic cohomology theory
was first introduced in mid 80's, see [Lan] and references therein,
though no direct
{\it geometric} interpretation of that theory is known as yet. Recently,
a definition of equivariant elliptic cohomology with complex
coefficients was suggested by Grojnowski [Gr].
His approach is not fully satisfactory however, since the construction
in [Gr] becomes void in the non-equivariant case: it cannot distinguish
non-equivariant elliptic cohomology from the ordinary cohomology at all.

In the present paper we adopt the following strategy.
In the first two sections we work out general formalism of
the `would-be' equivariant elliptic
cohomology theory as if such a theory existed.
Most of our attention is payed to those features of the theory that
are essentially different from the known cohomology theories, e.g.,
an invariant construction of the Chern classes. At the end
of $\S 1$ we sketch Grojnowski's definition (we slightly
modify his definition actually in order to remain in
the domain of algebraic, rather than complex-analytic, geometry).
Our exposition is thus very similar, in spirit, to the
Kazhdan-Lusztig exposition [KL], in the Hecke algebra
case.  There, basic properties of the equivariant
$K$-homology theory were listed in detail while the definition
itself was only sketched. The reason for this is the same both
in [KL] and in our paper:
the information about the cohomology theory under consideration
used in the study of elliptic, resp. Hecke, algebras
 does not depend in any way on the  definition
of the cohomology theory itself.

Our  new results are concentrated mostly in chapters 3-5.
The heart of the paper is the
algebro-geometric construction of n. 3.4, which is quite general
and does not depend on the rest of the material. The main results
of the parer are theorem 3.5, and theorem
5.13. Proof of proposition 3.6
is technically most complicated; it requires a refinement of
a lengthy argument used in [BLM].

Part 5 provides also
a self-contained exposition of elliptic algebras that
might be useful for the non-expert. The material involved
does not seem to be well covered in the literature.

The present paper should be regarded as a preliminary report on a wider
project. This project involves deeper understanding of the
recently discovered
relationship
between loop and `double-loop' groups on the one hand, and moduli
spaces of algebraic vector bundles on an elliptic curve and or an
elliptic
surface, on the other (cf. [EFK] and [GKV]). The subject is also
closely
related to the construction of `Springer resolutions' and theory
of `Springer representations' for loop groups.
These
matters will be discussed in more detail in a future publication.

We would like to thank I. Grojnowski,
M. Hopkins, Y. Manin, H. Miller and H. Nakajima for valuable
discussions and suggestions.

\vskip 17mm

\centerline {\bf 1. Equivariant elliptic cohomology : axiomatics.}

\vskip 1cm

\noindent {\bf (1.2) Elliptic curves.} Let $S$ be a scheme. By an elliptic
curve over $S$ we mean a group scheme $p: E\rightarrow S$ over $S$ whose
fiber at every geometric point is an elliptic curve. Thus $E$ is
equipped with the zero section $i: S\rightarrow E$. We denote the line
bundle $i^*\Omega^1_{E/S}$ on $S$ by $\omega_E$ or simply $\omega$.
We denote by $E^\vee \rightarrow S$ the dual elliptic curve
parametrizing  line bundles on $X$ of relative degree 0. Note that
$\omega_{E^\vee} = (\omega_E)^*$ and $E^{\vee\vee} = E$.
We denote by $E^r$ and $E^{(r)}$ respectively the Cartesian and symmetric
powers
of $E$ over $S$.

\vskip .3cm

\noindent {\bf (1.3) The space ${\cal X}_G$.}
Let $G$ be a compact Lie group. Then $G$ is a maximal compact
subgroup in a uniquely determined complex algebraic Lie group
$G_{\bf C}$, and
the categories of continuous finite-dimensional
representations of $G$ and algebraic
representations of $G_{\bf C}$ are known to be
equivalent. Moreover, it was shown
by Chevalley that $G_{\bf C}$ comes in fact from a group scheme
$G_{alg}$ over {\bf Z}.
For example, if $G = U(n)$ is the unitary group, then $G_{\bf C}
= GL_n({\bf C})$ and $G_{alg}$ is the group scheme $GL_n$.
If $G$ is a finite group, then $G_{alg} = G_{\bf C} = G$. We write
$G^0$ for the  identity component of  $G$, and the same
for  $G_{\bf C}$, $G_{alg}$.

Let $p: E\rightarrow S$ be an elliptic curve. A principal $G_{alg}$-bundle
$P$ over $E$ defines a finite covering  $\tilde E\rightarrow E$ consisting of
``connected components of the fibers". The canonical morphism $P\rightarrow
\tilde E$ makes $P$ a principal $G^0_{alg}$-bundle over $\tilde E$.

Let ${\cal X}_G$ be
the moduli scheme of semistable principal $G$-bundles $P$ over $E^\vee$
(see [R]) such that $P \rightarrow \tilde E$
is topologically trivial on each geometric fiber.
This is a smooth scheme over $S$, and we let $p_G: {\cal X}_G\rightarrow S$
denote the canonical projection. We write $0 \i {\cal X}_G$ for the
"zero section" given by the class of trivial $G$-bundle.
A homomorphism $\phi: G\rightarrow H$
of compact Lie groups induces a morphism of groups schemes
$G_{alg}\rightarrow H_{alg}$ and thus a morphism of
$S$-schemes ${\cal X}_\phi:
{\cal X}_G \rightarrow {\cal X}_H$.

We write$\omega$  for the line bundle $p_G^*\omega$ on ${\cal X}_G$,
for short.

\vskip .3cm

\noindent {\bf (1.4) Examples.}

\noindent {\bf (1.4.1)} If $G = S^1$ is the circle group,
then ${\cal X}_G= E^{\vee\vee} = E$.

\noindent {\bf (1.4.2)} If $G$ is a torus, then let
$\Gamma = {\rm Hom}( S^1, G)$ be the
lattice of 1-parameter subgroups in $G$, so $G = S^1 \otimes_{\bf Z}
\Gamma$. Then ${\cal X}_G = E\otimes_{\bf Z} \Gamma$ is the sum of
${\rm rk} (\Gamma)$ copies of $E$.

\noindent {\bf (1.4.3)} If $G = {\bf Z}/n$, then ${\cal X}_G = E_n \i E$
is the subgroup of $n$-torsion points of $E$.

\noindent {\bf (1.4.4)} The above examples can be unified as follows. Let $G$
be
a compact Abelian group. Then the Pontryagin dual
$\hat G = {\rm Hom}(G, S^1)$ is a discrete Abelian group, and
${\cal X}_{G} = {\rm Hom}(\hat G, E)$. In particular,
${\cal X}_G$ is naturally endowed with a structure of compact Abelian group.
Indeed, given an isomorphism of Lie groups $E\simeq S^1\times S^1$,
the Lie group ${\cal X}_G$ is naturally isomorphic to $G\times G$ in
such a way that for any closed subgroup $H\i G$ the closed subgroup
${\cal X}_H\i{\cal X}_G$ is identified with $H\times H\i G\times G$.

\noindent {\bf (1.4.5)} If $G$ is a connected compact Lie group with maximal
torus $T$ and Weyl group $W$, then ${\cal X}_G = {\cal X}_T/W$. For example, if
$G = U(n)$, then ${\cal X}_G = E^{(n)}$ is the $n$th symmetric product
of $X$ (over $S$). If $G = SU(n)$, then ${\cal X}_G$ is the kernel of
the addition map $E^{(n)} \buildrel + \over \rightarrow E$. This is
a ${\bf P}^{n-1}$-bundle over $S$. More generally, if $G$ is a simple Lie
group,
then, as shown by Bernstein and Schwartzman [BS], the fibers
of $p_G: {\cal X}_G \rightarrow S$ are weighted projective spaces.
The space ${\cal X}_G$ was also studied in [Lo].

\noindent {\bf (1.4.6)} Let $S= {\rm Spec} ({\bf C})$, so $E$ is an elliptic
curve over {\bf C}.
Then, by  Weil - Narasimhan - Seshadri theorem ${\cal X}_G$ is the same as
the moduli spaces of principal $G$-bundles over $E$ with an integrable
connection, i.e., of conjugacy classes of homomorphisms $\pi_1(E)
\rightarrow G$. Since $\pi_1(E) = {\bf Z}^2$, we find that ${\cal X}_G$ is
the set of conjugacy classes of pairs of commuting elements of $G$.
This is to be compared with [S2], [HKR].

\vskip .3cm

\noindent {\bf (1.5) Equivariant elliptic cohomology.}

For every scheme $Z$ we denote by ${\rm Coh}(Z)$ the category
of coherent  sheaves of ${\cal O}_Z$ - modules.

Let $G$ be a compact Lie group
By a $G$-cell we mean,  (see [W]),
a $G$-space of the form $(G/H)\times D^n$ where $H$ is a closed subgroup
in $G$ and $D^n$
is a disk of dimension $n$ with trivial $G$-action. A $G$-C{}W complex
(G-complex, for short)
is a $G$-space $M$ glued from $G$-cells. Thus
$M$ admits an increasing filtration $M_n$ by
$G$-subspaces such that $M_n/M_{n-1}$ is a bouquet of $G$-spaces of the
form $(G/H)\times S^n$, where $S^n$ is the $n$-sphere.
It is known that any smooth manifold with smooth $G$-action is a finite
$G$-complex, see [AP].

In this paper we choose a naive point of view on (equivariant)
cohomology theories as functors on pairs of
$G$-complexes, close to that of tom Dieck [Di].  Let $E\rightarrow S$
be an elliptic curve.
A $G$-equivariant elliptic cohomology theory
associated to $E$ is a collection of the following data :

\vskip .1cm

\noindent {\bf (1.5.1)} {\it Contravariant} functors
${\rm Ell}^i_G$ ($i=0,1,\ldots$) from the category of pairs of
finite $G$-complexes into ${\rm Coh}({\cal X}_G)$.
\vskip .1cm

\noindent {\bf (1.5.2)} {\it Natural} sheaf morphisms
$\partial: {\rm Ell}^i_G (A)
\rightarrow {\rm Ell}^{i+1}_G(M,A)$ given for any pair $(M,A)$.

\vskip .1cm

\noindent {\bf (1.5.3)} {\it Multiplication maps} ${\rm Ell}_G^i(M,A) \otimes
{\rm Ell}_G^j(N,B) \rightarrow {\rm Ell}_G^{i+j}(M\times N,
M\times B\cup A\times N)$ which are associative, graded commutative and
functorial in an obvious sense.

\vskip .2cm

\noindent Put ${\rm Ell}_G=\bigoplus_i {\rm Ell}^i_G$.
 The data (1.4.1-3) should satisfy
the following conditions :

\vskip .2cm

\noindent {\bf (1.5.4)} {\it Homotopy and exactness :} Two $G$-homotopic maps
of pairs induce the same maps on ${\rm Ell}^i_G$, and for every pair
$(M,A)$ the following natural sequence is exact
$$...\rightarrow {\rm Ell}^i_G(M,A)\rightarrow
{\rm Ell}^i_G(M)\rightarrow
{\rm Ell}^i_G(A)\buildrel \partial
\over\longrightarrow {\rm Ell}^{i+1}_G(M,A)
\rightarrow ...$$

\vskip .1cm

\noindent {\bf (1.5.5)} {\it Periodicity axiom :}
There are natural isomorphisms
$$
{\rm Ell}_G^{i-2}(M,A) \simeq {\rm Ell}^i_G(M,A)
\otimes \omega. \quad {\rm Moreover,}\quad
 {\rm Ell}_G^{2i+1}(pt) = 0, \quad {\rm Ell}_G^{2i}(pt) =
\omega^{\otimes (-i)}.$$

\vskip .1cm

Note that for the case $G=\{1\}$, $S=$ modular curve, (1.5.5) gives that the
space of global section of ${\rm Ell}_1^{-2i}(pt)$ is the space of
 modular forms
of weight $i$. Note also that (1.5.4) implies that ${\rm Ell}_G$ satisfies
the ``suspension axiom" : ${\rm Ell}^i_G(\Sigma M, \Sigma A) =
{\rm Ell}^{i-1}_G(M, A)$, where $\Sigma$ denotes suspension.

\vskip .3cm

\noindent {\bf (1.6) Relations between different groups.}
We now list axioms relating theories ${\rm Ell}_G$ for different $G$.

\vskip .1cm
\noindent{\bf (1.6.1)}
For any homomorphism of compact Lie groups $\phi: G\rightarrow H$,
there is a multiplicative morphism of functors
$T_\phi: {\rm Ell}_H \Rightarrow{{\cal X}_\phi}_*{\rm Ell}_G$
from the category of $H$-pairs of complex to ${\rm Coh}({\cal X}_H)$.
For any two composable homomorphisms $\phi,\psi$ we have
$T_{\psi\phi} = T_\phi\circ T_\psi$.

\vskip .1cm

\noindent {\bf (1.6.2)} {\it Change of groups :} If $\phi: G\rightarrow H$ is
a group homomorphism, then for any $H$-complex $M$ we have an Eilenberg-Moore
spectral sequence
$$E_1^{i,j} =
L_i{{\cal X}_\phi}^*{\rm Ell}^j_H(M) \Rightarrow {\rm Ell}^{j-i}_G(M).$$
Here $L_i{{\cal X}_\phi}^*$ is the $i$-th left derived functor of the inverse
image (This sequence converges, since ${\cal X}_G$ has finite Tor-dimension
over ${\cal X}_H$). In particular,

\vskip .1cm

\noindent {\bf (1.6.3)}
If ${\cal X}_\phi$ is a flat morphism, there is an isomorphism
${\cal X}_\phi^*{\rm Ell}_H\buildrel\sim\over\rightarrow {\rm Ell}_G.$

\vskip .1cm
\noindent{\bf (1.6.4)} {\it Induction axiom :} Let $H\i G$ be the embedding
of a closed normal subgroup, and $M$ be a $G$-complex such that the action
of $H$ is free. Denote by $p\,:\,M\rightarrow M/H$ and
$\phi\,:\,G\rightarrow G/H$ the projections. Then the map
$p^*\circ T_\phi(M/H)\,:\,
{\rm Ell}_{G/H}(M/H)\rightarrow{{\cal X}_\phi}_*{\rm Ell}_G(M)$
is an isomorphism of sheaves.

\vskip .1cm

\noindent {\bf (1.6.5)} {\it K\"unneth formula :} Let $G,H$ be two compact Lie
group, $M$ be a $G$-complex and $N$ be an $H$-complex. Then
${\rm Ell}_{G\times H} (M\times N) = {\cal X}_{\rho_G}^*{\rm Ell}_G(M)
\otimes {\cal X}_{\rho_H}^*{\rm Ell}_H(N)$, where $\rho_{_G}, \rho_{_H}$
are the projections of $G\times H$ to $G, H$.

\vskip .2cm
We conjecture that any elliptic curve $E$
gives rise to a unique equivariant
elliptic cohomology theory, natural in $E$. A more
systematic, if more technical, way of formulating this conjecture
is by using the language  of spectra representing cohomology theories.
Namely, for any $G$ there should exist a canonical sheaf
${\cal F}_G$ of $G$-spectra
(in the sense of [LMS]) on ${\cal X}_G$, and for any homomorphism
$\phi: G\rightarrow H$ we should have ${\cal F}_G = (L\phi^*){\cal F}_H$
where $L\phi^*$ is the derived functor of the inverse image.
For $G=\{1\}$ this would give a sheaf of spectra on $S$ considered in
[H]. In Section 2 we will give a construction valid over the rationals.

\vskip .3cm

\noindent {\bf (1.7) Some consequences of the definition.}\hfill\break
{\bf (1.7.1)}
For a locally compact space $M$ let $M\cup\infty$ denote the one-point
compactification of $M$. Then (1.5.4-5) implies that
(for the trivial $G$-action)
${\rm Ell}^0_G({\bf R}^n\cup\infty) = 0$ for odd $n$ while
${\rm Ell}_G^0({\bf C}^n\cup\infty) = \omega^{\otimes (n)}$.

\vskip .1cm

\noindent {\bf (1.7.2)} Let 1 denote the group with one element.
By applying the induction axiom (1.6.1) to finite-dimensional
$G$-invariant skeletons of the contractible space $EG$, we find that
${\rm Ell}_1^{0}(BG)$ is the completion of the sheaf
${\cal O}_{{\cal X}_G}$
at the point $0\in {\cal X}_G$. In particular, taking $G = S^1$,
we find that ${\rm Ell}_1^0 ({\bf CP}^\infty) := \lim \, {\rm inv}\,
{\rm Ell}_1^0 ({\bf CP}^n)$
is the completion of ${\cal O}_E$ at $0\in E$.

\vskip .1cm

\proclaim (1.7.3) Proposition. Let $\phi\,:\,H\hookrightarrow G$ be the
embedding of a closed subgroup of $G$. Then ${{\cal X}_\phi}_*{\rm Ell}_H(pt)$
and ${\rm Ell}_G (G/H)$ are isomorphic coherent sheaves on ${\cal X}_G$.

\noindent{\sl Proof.} Let $\delta\,:\,H\hookrightarrow G\times H$ be the
diagonal map and $\rho_{_G},\rho_{_H}$ be the projections of $G\times H$ to
$G$, $H$. Let
$\{1\}\buildrel i\over\rightarrow G\buildrel p\over\rightarrow\{1\}$
be the canonical maps. Since $p^*\circ T_{\rho_H}(pt)\,:\,
{\rm Ell}_H(pt)\rightarrow{{\cal X}_{\rho _H}}_*{\rm Ell}_{G\times H}(G)$
is an isomorphism by (1.6.4), the map $i^*\circ T_{\delta}(G)\,:\,
{\rm Ell}_{G\times H}(G)\rightarrow {{\cal X}_{\delta}}_*{\rm Ell}_H(pt)$
is either an isomorphism. Composing it with the isomorphism
${\rm Ell}_G(G/H)\rightarrow {{\cal X}_{\rho_G}}_*{\rm Ell}_{G\times H}(G)$
given by (1.6.4), one gets the result.

\vskip .1cm

\noindent {\bf (1.7.4) The scheme $M_{{\cal X}_G}$.}
By (1.5.3), any ${\rm Ell}^0_G(M)$ is a sheaf of commutative
algebras over ${\cal X}_G$. We denote its spectrum by $M_{{\cal X}_G}$
and denote by $\pi_{_M}: M_{{\cal X}_G}\rightarrow {\cal X}_G$
the projection. In particular,
${\rm Ell}^0_G(M)={\pi_{_M}}_*{\cal O}_{M_{{\cal X}_G}}$ and, if $M=\{ pt\}$,
we get $pt_{{\cal X}_G}={\cal X}_G$.
If $\phi\,:\,M\rightarrow N$ is a holomorphic map the contravariant
functoriality map $\phi^*\,:\,{\rm Ell}_G(N)\rightarrow{\rm Ell}_G(M)$
induces a morphism of schemes
$\phi_{{\cal X}_G}\,:\,M_{{\cal X}_G}\rightarrow N_{{\cal X}_G}$.
The maps $M\mapsto M_{{\cal X}_G}$, $\phi\mapsto\phi_{{\cal X}_G}$
define a functor from the category of $G$-complexes to the category
of schemes over ${\cal X}_G$. Note that, from
Proposition (1.7.3), if $M$ is a finite $G$-complex
then $\pi_{_M}$ is a finite morphism. Indeed Proposition (1.7.3)
can be rephrased in the following way :

\proclaim (1.7.5). Let $\phi\,:\,H\hookrightarrow G$ be the
embedding of a closed subgroup of $G$. Then $({\cal X}_H,{\cal X}_\phi)$
and $((G/H)_{{\cal X}_G},\pi_{_{G/H}})$ are isomorphic schemes over
${\cal X}_G$.

\noindent Section (1.6) can be rewritten in terms of $M_{{\cal X}_G}$.
Axioms (1.6.1) and (1.6.3) mean that a homomorphism
$\phi\,:\,G\rightarrow H$ and a $H$-complex $M$ give a commutative square
$$\matrix{&M_{{\cal X}_G}&\longrightarrow &M_{{\cal X}_H}&\cr
\pi_{_M}&\big\downarrow&&\big\downarrow&\pi_{_M}\cr
&{\cal X}_G&\buildrel {\cal X}_\phi\over\longrightarrow &{\cal X}_H&}$$
which is Cartesian if ${\cal X}_\phi$ is flat.
Similarly, (1.6.4) means that if a normal subgroup $H$ of $G$ acts freely
on a $G$-complex $M$ then the quotient map $\phi\,:\,G\rightarrow G/H$
induces a commutative diagram
$$\matrix{&M_{{\cal X}_G}&\buildrel =\over\longrightarrow &
(M/H)_{{\cal X}_{G/H}}&\cr
\pi_{_M}&\big\downarrow&&\big\downarrow&\pi_{_{M/H}}\cr
&{\cal X}_G&\buildrel {\cal X}_\phi\over\longrightarrow &{\cal X}_{G/H}&.}$$
Formally, $M_{{\cal X}_G}$, ${\cal X}_G$ and $\pi_{_M}$
behave exactly in the same way that $M_G$, $BG$ and the
canonical fibration $M_G\rightarrow BG$ do in a completed theory.

\vskip .1cm

\noindent {\bf (1.7.6)} For any triple $(M,A,B)$, $B\i A\i M$, of finite
$G$-complexes (1.5.4) and functoriality give an exact sequence
$$...\rightarrow {\rm Ell}^i_G(M,A)\rightarrow
{\rm Ell}^i_G(M,B)\rightarrow
{\rm Ell}^i_G(A,B)\buildrel \partial
\over\longrightarrow {\rm Ell}^{i+1}_G(M,A)
\rightarrow ...$$

\vskip .3cm

\noindent {\bf (1.8) Chern classes}.
If $H$ is another Lie group and $P\rightarrow M$ is a $G$-equivariant
principal $H$-bundle over a $G$-space $M$, then (1.6.1) gives a regular map
$$c_P: M_{{\cal X}_{G}} = P_{{\cal X}_{G\times H}}\rightarrow {\cal X}_{G\times
H}\rightarrow {\cal X}_H.$$
We call this map the {\sl characteristic class} of $P$. In particular, if
$V$ is a $G$-equivariant vector bundle over $M$ of rank $r$
with Hermitian form in the fibers, then the principal $U(r)$-bundle
of orthonormal frames in $V$ gives a map
$c_V: M_{{\cal X}_G}\rightarrow E^{(r)}$.
It is clear that
$$c_{V\oplus W} = c_{V}\oplus c_{W}, \quad
c_{V\otimes W} = c_{V}\otimes c_{W},$$
where $\oplus$ and $\otimes$ on the right hand sides of the equalities
are defined as follows :

$$\oplus\,:\,E^{(r_1)}\times E^{(r_2)} \rightarrow E^{(r_1+r_2)},\quad
{\rm the\quad symmetrization\quad map},$$

$$\otimes\,:\,E^{(r_1)}\times E^{(r_2)} \rightarrow E^{(r_1r_2)},\quad
(\{x_1,...,x_{r_1}\},\{y_1,...,y_{r_2}\})\mapsto \{x_i+y_j\}.$$

\vskip .2cm

By a {\it coordinate} on $E$ we mean a rational section
$f$ of the line bundle $p^*\omega^{-1}$ (where $p: E\rightarrow S$ is the
projection) with the following property:

\proclaim (1.8.1). The section $f$ is regular near $0 \i E$
and the differential $d_{E/S}(f)$, being restricted to
$0$, coincides with the identity map $\omega^{-1}\rightarrow
\omega^{-1}$.

Thus, on each geometric fiber $E_s, s\in S$, a coordinate $f$ gives a
rational function but with values in the tangent
space $T_0E_s$ (very much like the logarithm in a Lie group).
If $E$ is equipped with a coordinate, then for any $G$-equivariant
vector bundle $V$ on a $G$-space $M$ and any $i\geq 0$  we have its
equivariant Chern class
$c_i^f(V)$ which is a rational section of the  sheaf ${\rm Ell}_G^{2i}(M)$
regular near 0. The
class $c_i^f(V)$ is described as follows.
Namely, let $p_r: E^{(r)}\rightarrow S$ be the projection and
$\sigma_i(f)$ be the
rational section of the bundle $p_r^*\omega^{\otimes (-i)}$ on $E^{(r)}$
given by $\sigma_i(f)(x_1, ..., x_r) =\sigma_i(f(x_1),...,f(x_r))$,
where $\sigma_i$ is the $i$th elementary
symmetric function. We define $c_i^f$ to be the pullback of the
$\sigma_i(f)$
with respect to the map $c_V: M_{{\cal X}_G}\rightarrow E^{(r)}$.
This is a section of the bundle $\pi_M^*\omega^{\otimes (-i)}$
 on $M_{{\cal X}_G}$,
i.e., a section of ${\rm Ell}^{2i}_G(M)$ on ${\cal X}_G$.
It is clear that the Chern classes thus defined satisfy the usual
properties. In particular, it follows that the non-equivariant theory
${\rm Ell}_1$ becomes complex oriented. The corresponding formal group law
is the  series
$$F(x,y)=x+y+\sum_{i+j\geq 2} a_{ij}x^iy^j,
 \quad a_{ij}\in H^0(S, \omega^{\otimes i+j-1}) = {\rm Ell}_1^{-2i-2j+2}(pt),$$
giving the addition theorem
for  $f$, i.e. $f(u+v) = F(f(u), f(v))$.

A similar construction with $f$ a section of a non-trivial line bundle on $E$
gives the Euler class (see (2.6)).

\vskip .2cm

\noindent {\bf (1.8.2) Example.} Let $M=pt$ and $V={\bf C}^r$ be a vector
space,
$G=U(1)^n$ be a torus acting on $V$ via characters $\theta_1, ..., \theta_r$.
Let ${\cal X}_{{\theta}_i}\,:\,E^n\rightarrow E$ be the morphism
of algebraic groups given by the exponents entering into $\theta_i$. Then
$$c_V={\cal X}_{\theta_1}\oplus{\cal X}_{\theta_2}\oplus...
\oplus{\cal X}_{\theta_r}\,:\,E^n\rightarrow E^{(r)}.$$

\vskip .2cm

\noindent {\bf (1.8.3) Remarks.}

\noindent {\bf (1.8.3.1)} The formal germ at $0$ of a coordinate $f$
gives an element of ${\rm Ell}^2_1({\bf CP}^\infty)$ (just recall that
${\rm Ell}^0_1({\bf CP}^\infty)$ is the completion of ${\cal O}_E$
at $0\in E$) and the differential
$d_{E/S}f|_0$ along $0$  corresponds to the image of this element
in ${\rm Ell}_1^2({\bf CP}^1) = {\rm Ell}_1^0 (pt) = {\cal O}_S$.
So the above choice of the notion of coordinate exactly corresponds to
the concept of a complex orientation in topology, as explained in [AHS].

\noindent {\bf (1.8.3.2)} There are several versions of elliptic cohomology
considered in the literature, for instance, ``classical" elliptic cohomology
related to the Jacobi sine and points of order 2, Hirzebruch's level $N$
elliptic cohomology etc, see [Lan], [HBJ]. From our point of view,
the difference
among them is the artefact of the choice of coordinate,
which is made to pass from a formal group
( a group object in the category of formal schemes) to a formal group law
(a power series $F(x,y)$). Namely, an elliptic
curve  $E$ by itself does not have any preferred coordinate and so does not
define any formal group law. One natural choice of a coordinate is associated
to a point $\eta\in E$ of order 2. This is the function
(Jacobi sine) on $E$ having simple
zeroes at 0 and $\eta$ and simple poles at the two other points of
 second order. This choice of  coordinate leads to the level 2 elliptic
cohomology. In a similar way, a surjective homomorphism $\beta: E_N
\rightarrow {\bf Z}/N$, where $E_N$ is the set of $N$-torsion points on $E$,
gives a coordinate $f$ with simple zeroes on $\beta^{-1}(0)$ and simple
poles on $\beta^{-1}(1)$, and this choice of $f$ leads to the level N
elliptic cohomology of Hirzebruch [HBJ]. In our approach there is exactly one
elliptic cohomology theory associated to a given elliptic curve.

\vskip .3cm

\noindent {\bf (1.9) Cohomology of projective  and flag bundles.}
Let $M$ be a  finite $G$-complex and $V$ be a $G$-vector bundle on $M$,
 of rank $r$. We can assume that  $V$ is equipped with an Hermitian
metric in the fibers. Let ${\rm Fr}(V)$ be the space of orthonormal frames
in $V$. It is acted upon by $G\times U(r)$. The projectivization
${\bf P}(V)$ is the quotient
${\rm Fr}(V)/U(1)\times U(r-1)$. Thus (1.6.4) gives
$${\bf P}(V)_{{\cal X}_G} = {\rm Fr}(V)_{{\cal X}_{G\times U(1)\times U(r-1)}}
\quad , \quad M_{{\cal X}_G}={\rm Fr}(V)_{{\cal X}_{G\times U(r)}}.
\leqno {\bf (1.9.1)}$$
Consider the embedding $\phi: G\times U(1)\times U(r-1)\hookrightarrow
G\times U(r)$. From (1.7.5), the corresponding map
${\cal X}_\phi: {\cal X}_{G\times U(1)\times
U(r-1)}\rightarrow {\cal X}_{G\times U(r)}$ is just the symmetrization map
${\cal X}_G\times E\times E^{(r-1)}\rightarrow {\cal X}_G \times E^{(r)}$.
This map is flat and therefore
$${\bf P}(V)_{{\cal X}_G} =
{\rm Fr}(V)_{{\cal X}_{G\times U(1)\times U(r-1)}}=
{\rm Fr}(V)_{{\cal X}_{G\times U(r)}}
\times_{_{{\cal X}_{G\times U(r)}}}{\cal X}_{G\times U(1)\times U(r-1)}.
\leqno {\bf (1.9.2)}$$
Thus, since ${\rm Fr}(V)_{{\cal X}_{G\times U(r)}}=M_{{\cal X}_G}$,
we get the following Cartesian square
$$\matrix{{\bf P}(V)_{{\cal X}_G}&\longrightarrow& E\times E^{(r-1)}\cr
\big\downarrow&&\big\downarrow\cr
M_{{\cal X}_G}&\buildrel c_V\over\longrightarrow&E^{(r)}}.\leqno{\bf (1.9.3)}$$
Note that the map ${\bf P}(V)_{{\cal X}_G}\rightarrow M_{{\cal X}_G}\times E$
is
an embedding.

If we choose a coordinate $f$ on $E$ (by shrinking $S$, if necessary),
then we get Chern classes $c_i^f(V)\in \Gamma_{rat}({\rm Ell}^{2i}_G(M))$
and (1.9.3) yields, after passing to the completion at 0,
the standard description of the (completed) equivariant elliptic cohomology
of ${\bf P}(V)$ in terms of Chern classes.

In a similar way one can treat general flag varieties. Let ${\bf r} =
(r_1, ..., r_n)$ be a positive integer vector such that $\sum r_i=r$,
and let $F_{\bf r}(V)$ be the variety of flags of subspaces
 $V_1\i ... \i V_n$ in
fibers of $V$ such that $\dim(V_i/V_{i-1}) = r_i$. Let $E^{({\bf r})}=
\prod E^{(r_i)}$. Then we have a Cartesian square
$$\matrix{ F_{\bf r}(V)_{{\cal X}_G}&\longrightarrow &E^{({\bf r})}\cr
\big\downarrow&&\big\downarrow\cr
M_{{\cal X}_G}&\longrightarrow & E^{(r)}}, \leqno {\bf (1.9.4)}$$
where the right vertical map is the symmetrization.

\vskip .3cm

\noindent {\bf (1.10) Examples.}

\noindent {\bf (1.10.1)} Let $M=pt$ and $V={\bf C}^r$ be a vector space,
$G=U(1)^n$ be a torus acting on $V$ via
characters $\chi_1, ..., \chi_r$, so that $\chi_i$ is taken
$\mu_i$ times. Then ${\bf P}(V)_{{\cal X}_G} = \bigcup
\Gamma_i^{[\mu_i]}$ where  $\Gamma_i \i E^{n+1}$ is the graph of
the homomorphism of algebraic groups $\phi_i\,:\,E^n\rightarrow E$ given by
the exponents entering into $\chi_i$ and the superscript
$\mu_i$ means the $\mu_i$th infinitesimal neighborhood. The natural map
$${\bf P}(V)_{{\cal X}_G} = \bigcup
\Gamma_i^{[\mu_i]} \rightarrow E^n = {\cal X}_{G}$$
is just the natural projection of the graphs of the homomorphisms
onto the domain.

\vskip .2cm

\noindent {\bf (1.10.2)} Let $G=U(n)$, and $M$ be the projective space
${\bf CP}^{n-1}$ with standard action. Then
$\pi_{_M}\,:\,M_{{\cal X}_G}\rightarrow{\cal X}_G$ is the canonical projection
$E\times E^{(n-1)}\rightarrow E^{(n)}$.

\vskip .2cm

\noindent {\bf (1.10.3)}
Let $G=U(n)$ and $M$ be the $G$-manifold of complete flags in ${\bf C}^n$
with standard action. Then $\pi_{_M}\,:\,M_{{\cal X}_G}\rightarrow{\cal X}_G$
is just the canonical projection $E^n\rightarrow E^{(n)}$. Let $\theta$ be a
character of the torus $H=U(1)^n$ and $L_\theta=G\times_{_H}{\bf C}_\theta$ the
corresponding line bundle on $M$. From (1.7.3) and (1.8.2) the characteristic
class $c_{L_\theta}\,:\,E^n\rightarrow E$ is the homomorphism of algebraic
groups, ${\cal X}_\theta$, given by the exponents entering into $\theta$.

\vskip .3cm

\noindent {\bf (1.11) Grojnowski's construction.}
We sketch a construction of equivariant elliptic cohomology, ${\rm
Ell}^0_T(M)$,
for a compact torus $T= S^1\times\ldots\times S^1$ acting on a
a space $M$. In general, if $G$ is
 a connected compact group, $T$
a maximal  torus in $G$,
and
$M$ is a $G$-variety, we have by (1.4.5)
${\cal X}_G = {\cal X}_T/W$. Accordingly, we set by definition
$${\rm Ell}^0_G(M) := {\rm Ell}^0_T(M)^W,$$

Fix an elliptic curve $E$, viewed as a 1-dimensional complex
Lie group. Let ${\rm exp}: {\rm Lie} E \to E$ be the (non-algebraic)
exponential map.
A choice of basis in the lattice ${\rm Ker} {\rm exp}$
yields, via the exponential map, an isomorphism  $E\simeq\bf C/(\bf Z + \tau\bf
Z)$.
View ${\bf C} \simeq {\rm Lie} E$ as a
2-dimensional
real vector space $\bf R^2$. Then the isomorphism above
gives a real Lie group isomorphism
$$S^1\times S^1 \buildrel\sim\over\to E \quad,\quad
\bf R/\bf Z \times \bf R/{\bf Z }\ni (x, y) \mapsto x+
\tau\cdot y\; {\rm mod} \,(\bf Z + \tau\bf Z)
\eqno{\bf (1.11.1)}$$
The natural $SL_2(\bf Z)$-action  on  $\bf R^2$ induces an action on
$S^1\times S^1 \simeq \bf R/\bf Z \times \bf R/\bf Z$ by the formula:
$$
({a\atop c} {b\atop d}) : (u,v) \mapsto (u^a v^b, u^c v^d)\quad,\quad
u,v \in S^1$$
It is immediate from this formula
that the subgroup in $S^1$ generated by $u$ and $v$ depends only
on
the $SL_2(\bf Z)$-orbit of $(u,v)\in S^1\times S^1$ and not on a
particular representative. Hence, this subgroup is independent of the
choice
of basis in the  lattice ${\rm Ker} ({\rm exp})$. Thus, to any element
$e\in E$
we can associate canonically a subgroup in $S^1$ defined as follows
$$\langle e\rangle= \{u^a\cdot v^b\,,\,a,b\in \bf Z\},\eqno {\bf (1.11.2)})$$
 where the pair $(u,v)$ corresponds $e$ under isomorphism (1.11.1)

Now let $T$ be a compact torus, and
$\Gamma = {\rm Hom}( S^1, T)$ the
lattice of 1-parameter subgroups in $T$. Then, by $(1.4.2)$ we have
canonical isomorphisms
$$T= S^1 \otimes_{\bf Z}\Gamma\quad,\quad {\cal X}_T = E\otimes_{\bf
Z} \Gamma
\eqno {\bf (1.11.3)}$$
Combining these isomorphisms with
(1.11.1)
we get a chain of group  isomorphisms
$$E\otimes_{\bf Z} \Gamma \simeq (S^1\times S^1)\otimes_{\bf Z} \Gamma
\simeq T\times T \eqno {\bf (1.11.4)}$$
Let $x \in E\otimes_{\bf Z} \Gamma$, and let $(g_1,g_2)\in T\times T$
be the pair corresponding to $x$ via the above isomorphism.
The trivial, though crucial,
observation, based on formula (1.11.2),
is that,  the subgroup $\langle x\rangle\,\subset\,T$ generated by
$g_1$ and $g_2$
is canonically associated with $x$, and is independent of the
choice of basis in the kernel of ${\rm exp}: {\rm Lie} E \to E$.

Let $M$ be a $T$-space. For any  pair $g_1,g_2\in T$
write $M^{^{g_1,g_2}}$ for the set of points of $M$ simultaneously
fixed by $g_1$ and $g_2$. By the remark in the previous paragraph,
this fixed point set is totally  determined by the corresponding
point $x \in E\otimes_{\bf Z} \Gamma$, so that we may (and will) write
$M^{\langle x\rangle} := M^{^{g_1,g_2}}$

To define ${\cal X}_T(M)$ as a {\it complex analytic}
scheme we have, by (1.4.1),  to construct a
an {\it analytic} sheaf of algebras on $E\otimes_{\bf Z} \Gamma$ in the
Hausdorff topology.
The geometric fiber of this sheaf at a point $\langle x\rangle
\in E\otimes_{\bf Z} \Gamma$ will be
the finite dimensional $\bf C$-algebra $H^*(M^{\langle x\rangle})$,
the  ordinary
cohomology with complex coefficients of the simultaneous
fixed point set
of the pair $g_1,g_2\in T$ canonically associated to $x\in
{\rm Lie} E \otimes_{\bf Z} \Gamma$..
To glue the fibers together in an analytic sheaf,
suffices it to define, for any
point $x\in E\otimes_{\bf Z} \Gamma$
 and a small enough open neighborhood $U\,\subset\,E\otimes_{\bf Z} \Gamma$,
 the space $\Gamma(U,{\rm Ell}^0_T(M))$ of its holomorphic sections.

Let ${\bf t}_{_{\bf C}}$ be the complexified Lie algebra of the torus
$T$. By (1.11.3) we have a canonical isomorphism
$$({\bf t}_{_{\bf C}})
\simeq {\rm Lie} E \otimes_{\bf Z} \Gamma\eqno {\bf (1.11.5)}$$
Fix a point $x\in E\otimes_{\bf Z} \Gamma$ and its small connected neighborhood
$U$,
as in the previous paragraph. Let $\cal U$ be a connected component of the
inverse image
of $U$ under the natural exponential map
$$ {\rm exp}: {\bf t}_{_{\bf C}}=
({\rm Lie} E) \otimes_{\bf Z} \Gamma \to E \otimes_{\bf Z} \Gamma$$
This map gives an analytic isomorphism $\cal U
\buildrel{\rm exp}\over\longrightarrow U$, provided $U$ is small
enough.

Further, write $H^*_T(-)$ for the $T$-equivariant cohomology
functor. Recall that, for any $T$-space $Y$, the equivariant
cohomology $H^*_T(Y)$ has a natural $H^*_T(pt)$-module structure.
It is known that
$H^*_T(pt)=\bf C[{\bf t}_{_{\bf C}}]$, is the
polynomial algebra. We view the  polynomial algebra as a subalgebra
of the sheaf  $,{\cal O}^{an}_{{\bf t}_{_{\bf C}}}$, of holomorphic functions
on ${\bf t}_{_{\bf C}}$.

Now given a point $x\in E\otimes_{\bf Z} \Gamma$ and its small connected
neighborhood
$U$, following
Grojnowski, we set
$$\Gamma(U,{\rm Ell}^0_T(M)) :=\Gamma({\cal U} ,{\cal O}^{an}_{{\bf
t}_{_{\bf C}}})
\otimes_{_{\bf C[{\bf t}_{_{\bf C}}]}}
H^*_T(M^{\langle x\rangle}),$$
where $M^{\langle x\rangle} := M^{^{g_1,g_2}}$ is the simultaneous
fixed point set
of the pair $g_1,g_2\in T$ canonically associated to $x\in
({\rm Lie} E) \otimes_{\bf Z} \Gamma$. The arguments above show
that the RHS is well defined and is independent of the choices involved.

\vskip 2cm

\centerline {\bf 2. Thom bundles and Gysin maps.}

\vskip 1cm

\noindent {\bf (2.1) Thom bundles.}
Let $M$ be a finite $G$-complex and $V$ be a complex $G$-vector
bundle
over $M$ of rank $r$. The Thom space ${\rm Th}(V)$ is the quotient
${\bf P}(V\oplus {\bf C})/{\bf P}(V)$. The relative elliptic cohomology
${\rm Ell}_G^0({\bf P}(V\oplus {\bf C}),{\bf P}(V))$ is a sheaf of modules
over ${\rm Ell}_G^0(M)$ and thus has the form ${\pi_{_M}}_*\Theta(V)$ for
a uniquely determined coherent sheaf $\Theta(V)$ on $M_{{\cal X}_G}$.
This sheaf turns out to be invertible, and the following description of it
is immediately deduced from the long exact sequence of
cohomology and the description of cohomology
of projective bundles given in (1.9).

Consider ${\bf P}(V)_{{\cal X}_G}$ as a codimension one subscheme in
$M_{{\cal X}_G}\times E$. Let $D_V$
be the (effective) divisor in $M_{{\cal X}_G}\times E$ formed by irreducible
components of this subscheme with multiplicities naturally given by the scheme
structure. We denote by ${\cal O}(-D_V)$ the invertible sheaf on
$M_{{\cal X}_G}\times E$ whose sections are functions vanishing along $D_V$
(with the multiplicities prescribed by the divisor). Then
$$\Theta(V) = j^* {\cal O}(-D_V), \quad {\rm where}\quad
j:M_{{\cal X}_G}\rightarrow M_{{\cal X}_G}\times E\quad{\rm takes}\quad
x\mapsto
(x,0).\leqno {\bf (2.1.1)}$$
It follows that
$$\Theta(V\oplus V') = \Theta(V)\otimes\Theta(V'), \quad
\Theta(V^*) = \Theta(V). \leqno {\bf (2.1.2)}$$
The first of the above equalities allows us to extend $\Theta$
to a group homomorphism $K_G(M)\rightarrow {\rm Pic}(M_{{\cal X}_G})$, i.e.
to take $V$ to be a virtual vector bundle.
The product in cohomology (see (1.4.3)) and the contravariant functoriality
with respect to the projection ${\rm Th}(V)\rightarrow M$, see (1.4.1),
gives a sheaves morphism on ${\cal X}_G$
$${\rm Ell}^i_G(M)\otimes{\pi_{_M}}_* \Theta(V)\rightarrow{\rm Ell}_G^{i}({\rm
Th}(V)).
\leqno {\bf (2.1.3)}$$
One proves, glueing along the $G$-cells of $M$, that (2.1.3)
is an isomorphism.
It is the (twisted) ``Thom isomorphism".

In particular, when $V$ is trivial as a vector bundle, i.e. comes from
a representation $V$ of $G$, the Thom space ${\rm Th}(V)$ is
the ``$V$-suspension" $\Sigma^V(M)$ in the sense of tom Dieck [Di]
(see also [LMS]). So (2.1.3) replaces the datum (c) from Section 1 of
[Di].

\vskip .1cm

\noindent {\bf (2.1.4) Examples.}

\noindent {\bf (2.1.4.1)} Let $G=U(n)$ and $M$ be
the $G$-manifold of complete flags in ${\bf C}^n$
with standard action. Given a non-trivial character $\theta$
of the torus $H=U(1)^n$
let $L_\theta=G\times_{_H}{\bf C}_\theta$
be the corresponding line bundle on $M$ and
${\cal X}_\theta\,:\,E^n\rightarrow E$ be the homomorphism of algebraic groups
given by the exponents entering into $\theta$. Then (see (1.10.3))
$$\Theta(L_\theta)={\cal O}_{_{E^n}}(-{\rm Ker}\,{\cal X}_\theta)
\in{\rm Pic}\,(E^n).$$
Let
$M\buildrel i\over\hookrightarrow
{\rm Th}(L_\theta)\buildrel\pi\over\rightarrow M$
be respectively the inclusion of the zero section and the projection.
The pull-back morphism
$i^*\,:\,{\cal O}_{_{E^n}}(-{\rm Ker}\,{\cal X}_\theta)\rightarrow
{\cal O}_{_{E^n}}$ can be seen as a holomorphic section of
${\cal O}_{_{E^n}}({\rm Ker}\,{\cal X}_\theta)$.
Since $i^*$ admits a right inverse,
$\pi^*$, the divisor of this section is precisely ${\rm Ker}\,{\cal X}_\theta$.

\vskip.1cm

\noindent {\bf (2.1.4.2)} Let $G$ a compact Lie group and $M$
a connected manifold with trivial $G$-action. Thus,
$$\pi_{_M}\,:\,M_{_{{\cal X}_G}}=M_{_S}\times{\cal X}_{_G}\to{\cal X}_{_G}$$
is the second projection (here $S={\cal X}_{\{ 1\}}$ is the modular curve).
Let $V$ be a $G$-equivariant complex vector bundle of rank $r$
on $M$. Since the character group of $G$ is discrete,
the action of $G$ on the fibers of $V$ is constant.
Denote ${\bf C}^r_{_\theta}$ the representation of $G$
on a fiber of $V$, and $\theta\,:\,G\to U(r)$
the corresponding morphism of groups.
Let $\Theta({\bf C}^r_{_\theta})\in{\rm Pic}({\cal X}_G)$
be the Thom bundle. Then, $\Theta(V)=\pi_{_M}^*\Theta({\bf C}^r_{_\theta})$.

\vskip .3cm

\noindent {\bf (2.2) Tubular neighborhood and elliptic homology.}
Let $(M,N)$ be a finite $G$-pair. We denote by $\Theta^i_N(M)$
the coherent sheaf on $N_{{\cal X}_G}$ such that ${\pi_{_N}}_*
\Theta^i_N(M) = {\rm Ell}^i_G(M, M\setminus N)$.
(We can replace $N$ by a small $G$-stable tubular neighborhood $T$ and
consider ${\rm Ell}^i_G(M, M\setminus {\rm Int}(T))$,
so as not to leave the category of CW-pairs).

In the particular case when $N\hookrightarrow M$ is an embedding
of complex manifolds with $G$ acting by holomorphic transformations,
$\Theta^{2i}_N(M) =\Theta(T_NM)\otimes\omega^{\otimes (-i)} $ is the
twisted Thom sheaf of the normal bundle of $N$ in
$M$, and $\Theta_N^{2i+1}(M) = 0$. In general, the graded sheaf
$\Theta_N(M)=\bigoplus_i\Theta_N^i(M)$ is the elliptic
analog of the Borel-Moore homology of $N$. It depends, however, on
the ambient space $M$. In particular,
if $M$ is a $G$-manifold and $F\i X\i M$ are closed
$G$-subcomplexes of $M$, the analogue of the long exact sequence
in homology associated to $X$, $F$ and $U=X\setminus F$ is the following
exact sequence of sheaves on ${\cal X}_G$ (see (1.7.6)) :
$$...\rightarrow
{\pi_{_F}}_*\Theta_F^i(M)\rightarrow
{\pi_{_X}}_*\Theta_X^i(M)\rightarrow
{\pi_{_U}}_*\Theta_U^i(M\setminus F)\buildrel \partial
\over\longrightarrow
{\pi_{_F}}_*\Theta_F^{i+1}(M)\rightarrow ...
\leqno{\bf (2.2.1)}$$

We have an obvious map of sheaves on ${\cal X}_G$
$${\pi_{_N}}_*\Theta_N(M)={\rm Ell}_G(M,M\setminus N)\rightarrow{\rm
Ell}_G(M)$$
which  (for smooth $N$) can be seen as a particular instance of direct
functoriality (Gysin map) for elliptic cohomology.
Let us consider this in general.

\vskip .3cm

\noindent {\bf (2.3) Gysin maps.}
Let $\phi: M\rightarrow N$ be a holomorphic map of complex
$G$-manifolds, and $d = \dim\, M - \dim\, N$ be the relative dimension
of $\phi$. The contravariant functoriality with respect to $\phi$
gives a morphism of schemes
$\phi_{_{{\cal X}_G}}: M_{{\cal X}_G}\rightarrow N_{{\cal X}_G}$
(see (1.7.4)). As such, it consists of
a morphism of topological spaces (still denoted $\phi_{_{{\cal X}_G}}$)
and a morphism
$\phi_{_{{\cal X}_G}}^{-1}{\cal O}_{N_{{\cal X}_G}}\rightarrow
{\cal O}_{M_{{\cal X}_G}}$
of sheaves on $M_{{\cal X}_G}$.
Let us denote by $\Theta(\phi)$ the Thom sheaf on $M_{{\cal X}_G}$
corresponding to the virtual vector bundle $\phi^*TN - TM$ on $M$.
If $\phi$ is proper, the direct image (or Gysin map) is a morphism $$\phi_*:
\Theta(\phi)\rightarrow \phi_{_{{\cal X}_G}}^{-1}{\cal O}_{N_{{\cal X}_G}}
\leqno {\bf (2.3.1)}$$
of sheaves on $M_{{\cal X}_G}$.
The construction is achieved in a standard way, by factoring $\phi$
through $G$-equivariant complex oriented embedding into a Euclidean space
and then considering the tubular neighborhoods [LMS].
Equivalently, we can view $\phi_*$ as a morphism
$$\phi_*: {\phi_{_{{\cal X}_G}}}_* \Theta(\phi)\rightarrow {\cal O}_{N_{{\cal
X}_G}}
\leqno {\bf (2.3.2)}$$
of sheaves on $N_{{\cal X}_G}$. This last map of sheaves is a homomorphism of
${\cal O}_{N_{{\cal X}_G}}$ -modules. This condition is the analog of
the projection
formula for the ordinary Gysin morphism.
Taking the projection to ${\cal X}_G$, one gets the following morphism
of sheaves on ${\cal X}_G$
$${\pi_{_M}}_*\Theta(\phi)\rightarrow{\pi_{_N}}_*{\cal O}_{N_{{\cal X}_G}} =
{\rm Ell}_{G}(N)\leqno {\bf (2.3.3)}$$

\vskip .3cm

\noindent {\bf (2.4) The direct image of a composition.}
The formula $(\psi\phi)_* = \psi_*\phi_*$ for direct images in
traditional theory is modified in our sheaf theoretic approach as follows.
Let $M\buildrel\phi\over\rightarrow N \buildrel
\psi\over\rightarrow L$
be two morphisms of complex $G$-manifolds. Then
$$\Theta(\psi\phi) = \Theta(\phi)\otimes
\phi_{_{{\cal X}_G}}^*\Theta(\psi),\leqno{\bf (2.4.1)}$$
and we have a natural commutative diagram :
$$\matrix{& {\pi_{_M}}_*(\Theta(\phi)\otimes
\phi_{_{{\cal X}_G}}^*\Theta(\psi)) &\buildrel =\over\longrightarrow &
{\pi_{_M}}_*\Theta(\psi\phi)&\cr
\phi_*\otimes{\rm Id}&\big\downarrow&&\big\downarrow&(\psi\phi)_*\cr
&{\pi_{_N}}_*\Theta(\psi)&\buildrel\psi_*\over\longrightarrow&
{\rm Ell}_G(L)&}\leqno{\bf (2.4.2)}$$

\noindent A little bit more transparent view of direct images and their
composition properties can be obtained by a certain twisting described
in the following

\proclaim (2.4.3) Proposition.
For every commutative diagram of proper morphisms of complex $G$-manifolds
$$\matrix{&M&\buildrel\phi\over\longrightarrow &N&\cr
\alpha&\quad\searrow&&\swarrow\quad&\beta\cr
&&Y&&}$$ the direct image map $\phi_*$ gives rise to a morphism
$\phi_*: {\pi_{_M}}_*\Theta(\alpha)\longrightarrow
{\pi_{_N}}_*\Theta(\beta)$. The correspondence $M\mapsto
{\pi_{_M}}_*\Theta(\alpha)$ extends to a covariant functor
from the category of complex $G$-manifolds over $Y$ (in which morphisms
are commutative triangles) to the category of sheaves on ${\cal X}_G$.

\vskip .1cm

\noindent {\bf (2.4.4) Remark.} Let $f\,:\,M\to N$ be a map between to
$G$-manifolds $M$ and $N$ and let $V$ be a $G$-vector bundle on $N$.
Denote by $i\,:\,N\hookrightarrow V$ the zero section. By applying
(2.4.1) to $M\buildrel f\over\to N\buildrel i\over\to V$ we get
$$f^*\Theta(V)=f^*\Theta(i)=\Theta(i\,f)\otimes\Theta(f)^{-1}.$$
Then, a simple computation gives : $f^*\Theta(V)=\Theta (f^*V)$.
In other words, $\Theta$ is a morphism of contravariant functors.

\vskip .3cm

\noindent {\bf (2.5) The singular case.}
We define Thom sheaf and Gysin map for any  proper map
possibly singular complex algebraic $G$-varieties by including
$\psi$ into a commutative diagram
$$ \matrix{& S&\buildrel\psi\over\rightarrow &T&\cr j&\big\downarrow &&
\big\downarrow &k\cr
&M&\buildrel \phi\over \rightarrow &N&}\leqno {\bf (2.5.1)}$$
where $j,k$ are closed embeddings, $M,N$ are smooth. Then we define
$$\Theta(\psi) = \Theta(j)\otimes
{\cal H}om\bigl(\psi_{_{{\cal X}_G}}^* \Theta(k),j_{_{{\cal X}_G}}^*
\Theta(\phi)\bigl).\leqno {\bf (2.5.2)}$$
It is easy to verify that such a definition of $\Theta(\psi)$ is indeed
independent of the choice of $j,\phi$ and $k$.
All the statements of (2.2-4) extend to this case
without any  further modifications. Note that with this definition,
for {\rm any} diagram of proper morphisms of form (2.5.1)
we will have equality (2.5.2).

\vskip .3cm

\noindent{\bf (2.6) Euler class.} Let $M$ be a finite $G$-complex and $V$ be
a complex $G$-vector bundle over $M$. Let $i:M\hookrightarrow V$ be the zero
section. The morphism of schemes
$i_{_{{\cal X}_G}}: M_{{\cal X}_G}\rightarrow V_{{\cal X}_G}$
is an isomorphism by the homotopy invariance of elliptic cohomology.
The map $i$ is proper. The Gysin map
$i_*\,:\,
\Theta(V)\rightarrow i^{^{-1}}_{_{{\cal X}_G}}{\cal O}_{V_{{\cal X}_G}}=
{\cal O}_{M_{{\cal X}_G}}$
can be seen as an regular section of the line bundle
$\Theta(V)^{^{-1}}$ on $M_{{\cal X}_G}$. Denote this section by
$e_V\in H^0(M_{{\cal X}_G},\Theta(V)^{^{-1}})$. In our situation, $e_V$ is
the analogue of the usual Euler class in equivariant cohomology or K-theory.
In particular, if $V$ and $W$ are complex vector bundles on $M$ one has
$e_{V\oplus W}=e_V\otimes e_W$ as a section of
$\Theta(V\oplus W)^{^{-1}}=\Theta(V)^{^{-1}}\otimes\Theta(W)^{^{-1}}$.

\vskip.1cm

\noindent{\bf (2.6.1) Example.} Let $\theta\,:\,G\to U(1)$ be a
character of $G$. Take $M=pt$ and $V={\bf C}_{_\theta}$ the corresponding
one-dimensional representation of $G$. The characteristic class of
${\bf C}_{_\theta}$ is the group homomorphism
${\cal X}_{_\theta}\,:\,{\cal X}_{_G}\to{\cal X}_{_{U(1)}}=E$.
Thus, $\Theta({\bf C}_{_\theta})={\cal O}(-{\rm Ker}\,{\cal X}_{_\theta})
\in{\rm Pic}({\cal X}_{_G})$.
The same argument as in (2.1.4.1) shows that the Euler class
$e_{_{{\bf C}_\theta}}$ is a regular section of
${\cal O}({\rm Ker}\,{\cal X}_{_\theta})$ over ${\cal X}_{_G}$
with zeros precisely on ${\rm Ker}\,{\cal X}_{_\theta}$. Note that
${\rm Ker}\,{\cal X}_{_\theta}={\cal X}_{_{{\rm Ker}\,\theta}}$
(see (1.4.4)).

\vskip .3cm

\noindent {\bf (2.7) Action of correspondences on elliptic cohomology.}
By a correspondence between smooth complex algebraic
$G$-varieties $M_1$ and $M_2$ we
mean a $G$-stable (possibly singular)
subvariety $W\i M_1\times M_2$ such that the
projection $q_2: W \rightarrow M_2$ is proper.
Any correspondence defines a morphism of sheaves on ${\cal X}_G$
$${\pi_{_W}}_*\Theta(q_2)\longrightarrow
{\cal H}om_{_{{\cal X}_G}}({\rm Ell}_G(M_1),\, {\rm Ell}_G(M_2)),
\leqno {\bf (2.7.1)}$$
which is defined as follows. Given any (local) section $s$ of
$\Theta(q_2)$ on $W_{{\cal X}_G}$ and a function $g$ on ${M_1}_{{\cal X}_G}$
(i.e., a section of ${\rm Ell}_G(M_1)$), we first lift $g$ to a function on
$W_{{\cal X}_G}$ by means of the scheme map
${q_1}_{_{{\cal X}_G}}: W_{{\cal X}_G}\rightarrow{M_1}_{{\cal X}_G}$
associated to the first projection
$q_1\,:\, W \to M_1$. Multiplying by $s$,
yields another section of
$\Theta(q_2)$, and applying ${q_2}_*$ we get a function
on ${M_1}_{{\cal X}_G}$.

\vskip .2cm

Given two correspondences $W_{12}\i M_1\times M_2$ and $W_{23}\i M_2\times
M_3$,
define their composition as the correspondence
$$W_{12}\circ W_{23} = \bigl\{(m_1,m_3)\in M_1\times M_3 \;|\;\,
\exists m_2\in M_2 \,\, {\rm s.t.} \, (m_1,m_2)\in W_{12},\,
(m_2,m_3)\in W_{23}\bigl\}.$$
Let $q_2: W_{12} \to M_2$, $q_3: W_{23} \to M_3$ and
$\overline q_3: W_{13}=W_{12}\circ W_{23} \to M_3$ be the second projections.
The composition of actions of correspondences on ${\rm Ell}_G$
lifts to a morphism of sheaves on ${\cal X}_G$
$${\pi_{_{W_{12}}}}_*\Theta(q_{2})\otimes
{\pi_{_{W_{23}}}}_*\Theta(q_{3})\rightarrow
{\pi_{_{W_{13}}}}_* \Theta(\overline q_{3})\leqno {\bf (2.7.2)}$$
which is defined as follows. Consider first $W_{12}\times_{_{M_2}}W_{23}$, the
fiber product of $W_{23}$ and $W_{12}$ over $M_2$. Then we have the diagram

$$\matrix{& W_{12}\times_{_{M_2}}W_{23}&\buildrel q_{23}\over\longrightarrow&
W_{23}&\buildrel q_{3}\over\longrightarrow M_3\cr
q_{12}&\big\downarrow&&\big\downarrow&&\cr
&W_{12}&\buildrel q_2\over\longrightarrow&M_2&&\cr
q_1&\big\downarrow&&&&\cr
&M_1&&&&}$$
The diagram yields
$$\Theta(q_{3}q_{23}) = q_{12}^*\Theta(q_{2})\otimes q_{23}^*\Theta(q_{3}),$$
and thus we have a multiplication
$${\pi_{_{W_{12}}}}_*\Theta(q_{2})\otimes {\pi_{_{W_{23}}}}_*\Theta(q_{3})
\rightarrow {\pi_{_{W_{12}\times_{_{M_2}}W_{23}}}}_* \Theta(q_{3}q_{23}).$$
The map (2.7.2) is obtained by composing this multiplication with the
map
$${\pi_{_{W_{12}\times_{_{M_2}}W_{23}}}}_* \Theta(q_{3}q_{23})\rightarrow
{\pi_{_{W_{13}}}}_* \Theta(\overline q_{3})$$
induced, as in (2.4.3), by the natural projection $W_{12}\times_{_{M_2}}W_{23}
\rightarrow W_{13}$.
The maps (2.7.2) are associative in an obvious sense. In particular,
if $M_1=M_2=M_3=M$ we get

\proclaim (2.7.3) Proposition. If $W\i M\times M$ is a correspondence such
that $W\circ W=W$, then ${\pi_{_W}}_*\Theta(q_2)$ has a natural
structure of a sheaf of associative algebras on ${\cal X}_G$.

\vskip.1cm

\noindent{\bf (2.7.4) Remark.} The morphism (2.7.1) is nothing but the
particular case of (2.7.2) corresponding to $M_1=pt$.

\vskip .3cm

\noindent{\bf (2.8) Lagrangian correspondences and their action on
${\rm Ell}_G$.} We will apply
a slight modification of the formalism above to the cotangent bundles.
Given a complex manifold $M$, let $TM$ (resp. $T^*M$) denote its
tangent (resp. cotangent) bundle. Write
$T_WM$ and $T^*_WM$ for the normal and conormal bundles of a
submanifold
$W\i M$. It is well
known that $T^*M$ is a symplectic manifold, and $T^*_WM$ is a
Lagrangian submanifold in $T^*M$.

Let $M_1,M_2$ be two smooth complex algebraic $G$-varieties
$N_1=T^*M_1$, $N_2=T^*M_2$.
By a Lagrangian
correspondence between $N_1$ and $N_2$ we mean a
$G$-stable possibly singular Lagrangian subvariety
$Z\i T^*(M_1\times M_2)$ whose projections to
$T^*M_1, T^*M_2$ are proper. In this case we have the diagram
$$M_1\buildrel\pi_1\over\longleftarrow T^*M_1\buildrel p_1\over
\longleftarrow Z\buildrel p_2\over\longrightarrow
T^*M_2\buildrel i_2\over\longleftarrow M_2,\leqno {\bf (2.8.1)}$$
where $\pi_1$ is the projection of the cotangent bundle, and $i_2$ is
the zero section. To a Lagrangian correspondence $Z$ we associate the
coherent sheaf
$$\Xi_{_Z} = \Theta(p_2)\otimes
\Theta(p_1^*\pi_1^*T^*M_1)\otimes\Theta(p_2^*\pi_2^*T^*M_2)^{-1}=
\leqno{\bf (2.8.2)}$$
$$=\Theta_{_Z}(T^*(M_1\times M_2))\otimes\Theta(p_1^*\pi_1^*T^*M_1)^{-1}\otimes
\Theta(p_2^*\pi_2^*T^*M_2)^{-1}$$
on $Z_{{\cal X}_G}$ called the {\it microlocal Thom sheaf} of $Z$. There is a
map
$${\pi_{_Z}}_*\Xi_{_Z}\rightarrow {\cal H}om_{_{{\cal X}_G}} ({\pi_{_{M_1}}}_*
\Theta(T^*M_1)^{-1}, {\pi_{_{M_2}}}_*\Theta(T^*M_2)^{-1})
\leqno {\bf (2.8.3)}$$
defined as follows. Suppose we have a local section of
$\Xi_{_Z}$ which we take in the form $a\otimes b \otimes c$
according to the factorization in the first line of (2.8.2).
Given a section $s$ of $\Theta(T^*M_1)^{-1}$,
the product $sb$ is a function on ${M_1}_{_{{\cal X}_G}}$,
we lift it to $Z_{_{{\cal X}_G}}$ by ${\pi_1}_{_{{\cal X}_G}}$
and ${p_1}_{_{{\cal X}_G}}$, then multiply by $a$,
getting a section of $\Theta(p_2)$, project along $p_2$,
getting a function on $(T^*M_2)_{_{{\cal X}_G}}$,
pull it back to ${M_2}_{_{{\cal X}_G}}$ by ${i_2}_{_{{\cal X}_G}}\,:\,
{M_2}_{_{{\cal X}_G}}\buildrel\sim\over\longrightarrow
(T^*M_2)_{_{{\cal X}_G}}$,
and multiply it with $c$, getting a section of $\Theta(T^*M_2)^{-1}$.

As for the non-modified case (2.7.2), given two Lagrangian correspondences
$Z_{12}\i N_1\times N_2$, $Z_{23}\i N_2\times N_3$, and their
composition $Z_{13}=Z_{12}\circ Z_{23}$, the composition of actions
of correspondences on Thom bundles as in (2.8.3) lifts
to an associative morphism of sheaves
$${\pi_{_{Z_{12}}}}_*\Xi_{_{Z_{12}}}\otimes{\pi_{_{Z_{23}}}}_*\Xi_{_{Z_{23}}}
\rightarrow{\pi_{_{Z_{13}}}}_*\Xi_{_{Z_{13}}}.\leqno{\bf (2.8.4)}$$
In particular, if $M_1=M_2=M$ we get the following :

\proclaim (2.8.5) Proposition. If $Z\i T^*M\times T^*M$ is a Lagrangian
correspondence such that $Z\circ Z=Z$, then ${\pi_{_Z}}_*\Xi_{_Z}$
has a natural structure of a sheaf of associative algebras on ${\cal X}_G$,
and ${\pi_{_M}}_*\Theta(T^*M)^{-1}$ is a sheaf of its modules.

A natural example of a Lagrangian correspondence is obtained as follows.
Let $W\i M_1\times M_2$ be a smooth closed submanifold whose projections
$q_1$ and $q_2$ to $M_1$ and $M_2$ are proper.
Take $Z = T^*_W(M_1\times M_2)$. Then $Z$ is
a Lagrangian correspondence between $N_1$ and $N_2$ while $W$ is a
correspondence between $M_1$ and $M_2$. We have

\proclaim (2.8.6) Proposition. Under the above assumptions, the sheaf
$\Xi_{_Z}$ is naturally identified with
${\cal O}_{Z_{{\cal X}_G}}$.

\noindent {\sl Proof.} The zero section map
$i\,:\, W\rightarrow Z = T^*_W(M_1\times M_2)$
is a $G$-homotopy equivalence, so
$i_{_{{\cal X}_G}}\,:\,W_{{\cal X}_G}\rightarrow Z_{{\cal X}_G}$
is an isomorphism of schemes. Thus it suffices to proove that
$i_{_{{\cal X}_G}}^*\Theta(p_2)\simeq {\cal O}_{W_{{\cal X}_G}}$.
Recall that $\Theta$ factors through
$K_G(W)$, the Grothendieck group of $G$-vector bundles on $W$.
Denoting by $[V]$ the class in $K_G(M)$ of a bundle $V$, we find
$$q_1^*[T^*M_1]=[T^*q_2]+[T^*_{_W}(M_1\times M_2)].$$
Since $\Theta$ is insensitive to the dualization (see (2.1.2)) we get
$$\Theta(q_2)\otimes\Theta(q_1^*T^*M_1)\otimes\Theta(i)^{-1}=
{\cal O}_{W_{{\cal X}_G}}.\leqno{\bf (2.8.7)}$$
In another hand, formula (2.4.1) applied to the commutative square
$$\matrix{&Z&\buildrel p_2\over\rightarrow&N_2&\cr
i&\uparrow&&\uparrow&i_2\cr
&W&\buildrel q_2\over\rightarrow&M_2&}$$
gives $\Theta(q_2)\otimes{q_2}_{_{{\cal X}_G}}^*\Theta(i_2)=
\Theta(i)\otimes i_{_{{\cal X}_G}}^*\Theta(p_2)$.
Combining it with (2.8.2) and (2.8.7) we get
$i_{_{{\cal X}_G}}^*\Xi_{_Z} = {\cal O}_{W_{{\cal X}_G}}$.

\vskip .3cm

\noindent{\bf (2.8.8) Remarks.}

\noindent{\bf (2.8.8.1)} If $W\i M_1\times M_2$ is a smooth locally closed
submanifold and $Z=T^*_W(M_1\times M_2)$ set
$\partial Z=\overline Z\setminus Z$ and define
$$\Xi_Z=\Theta_{_Z}(N_1\times N_2\setminus\partial Z)\otimes\Theta(p_1^*\pi_1^*
T^*M_1)^{-1}\otimes \Theta(p_2^*\pi_2^*T^*M_2)^{-1}.$$
Then we still have $\Xi_Z={\cal O}_{Z_{{\cal X}_G}}$.

\noindent{\bf (2.8.8.2)} Let $Z'$ be a $G$-stable Lagrangian
subvariety of $Z$. Let $i\,:\,Z'\rightarrow Z$ and $p'_2\,:\,Z'\rightarrow M_2$
be the inclusion in $Z$ and the projection to $M_2$. Formula (2.4.1) gives
$$\Xi_{_{Z'}}=\Theta(i)\otimes i^*_{_{{\cal X}_G}}\Xi_{_Z}.$$
Thus, the Gysin map induces a morphism
${i_{_{{\cal X}_G}}}_*\Xi_{_{Z'}}\rightarrow\Xi_{_Z}$
which commutes with the action of ${\pi_{_{Z}}}_*\Xi_{_Z}$ and
${\pi_{_{Z'}}}_*\Xi_{_Z'}$ on elliptic cohomology, as described in (2.8.3).

\noindent{\bf (2.8.8.3)} The morphism (2.8.3) is nothing but the particular
case of (2.8.4) corresponding to $M_1=N_1=pt$.

\vskip .3cm

\noindent{\bf (2.9) Localization.} Suppose that $G$ is a compact Abelian
Lie group and $M$ a complex smooth manifold. For any $g\in G$ put
${\cal X}_g={\cal X}_G\setminus\bigcup_{_{g\notin H}}{\cal X}_H$,
the union beeing taken over all the closed subgroups $H$ of $G$ not containing
$g$; the set ${\cal X}_g$ is non-empty (see example (1.4.4)).
The fixed-points set $M^g$ is a smooth complex subvariety of $M$.
Let $i\,:\,M^g\hookrightarrow M$ be the inclusion map; it is a proper map.

\proclaim (2.9.1) Proposition. For any $g\in G$, the Gysin map
$i_*\,:\,\pi_{_{M^g *}}\Theta(T_{_{M^g}}M)\rightarrow {\rm Ell}_G(M)$
is an isomorphism over ${\cal X}_g\i{\cal X}_G$.
Besides, the morphism of sheaves
$i^*i_*\,:\,\pi_{_{M^g *}}\Theta(T_{_{M^g}}M)\to {\rm Ell}_G(M^g)$
is the multiplication by the Euler class $e_{_{T_{M^g}M}}$ and is
invertible over ${\cal X}_g\i{\cal X}_G$.

\noindent{\sl Proof}. Denote by $j\,:\,M\setminus M^g\rightarrow M$ the
inclusion. The long exact sequence (2.2.1) gives :
$$...\buildrel\partial\over\rightarrow\pi_{_{M^g *}}\Theta^i(T_{_{M^g}}M)
\buildrel i_*\over\rightarrow {\rm Ell}^i_G(M)
\buildrel j^*\over\rightarrow {\rm Ell}^i_G(M\setminus M^g)
\buildrel\partial\over\rightarrow\pi_{_{M^g *}}\Theta^{i+1}(T_{_{M^g}}M)
\rightarrow ...$$
So it suffices to prove that, as a sheaf over ${\cal X}_G$,
${\rm Ell}^0_G(M\setminus M^g)$ is supported on ${\cal X}_G\setminus{\cal
X}_g$.
By a standard argument we are reduced to the case where $M\setminus M^g$
is a $G$-orbit, say $G/H$, with $g\notin H$. Then, we know from the induction
axiom (1.6.2) that ${\rm Ell}^0_G(G/H)={\cal O}_{{\cal X}_H}$.
Thus, the first part of the proposition follows from
${\cal X}_H\cap{\cal X}_g=\emptyset$. The second part of the proposition
is a direct corollary of the construction of Euler classes (see (2.6)),
except for the invertiblity over ${\cal X}_g$ which follows from standard
arguments and (2.6.1).

\vskip.2cm

Let $f\,:\,M_1\to M_2$ be a proper morphism of complex smooth varieties.
Proposition (2.4.3) applied to the diagram

$$\matrix{&M^g_1&\buildrel i_1\over\hookrightarrow&M_1&\cr
f^g&\downarrow&&\downarrow&f\cr
&M^g_2&\buildrel i_2\over\hookrightarrow&M_2&},$$
where $i_{_1}$, $i_{_2}$ are the inclusions and $f^g$ is the restriction of
$f$ to $M^g_1$ and $M^g_2$, gives the following commutative square
$$\matrix{{\pi_{_{M^g_1}}}_*\Theta(TM^g_1)^{^{-1}}&\buildrel{i_1}_*\over
\longrightarrow&{\pi_{_{M_1}}}_*\Theta(TM_1)^{^{-1}}\cr
f^g_*\downarrow&&\downarrow f_*\cr
{\pi_{_{M^g_2}}}_*\Theta(TM^g_2)^{^{-1}}&\buildrel{i_2}_*\over\longrightarrow
&{\pi_{_{M_2}}}_*\Theta(TM_2)^{^{-1}}}\leqno{\bf (2.9.2)}$$
where ${i_{_1}}_*$ and ${i_{_2}}_*$ are given by multiplication by
$e_{_{T_{M^g_1}M_1}}$ and $e_{_{T_{M^g_2}M_2}}$. The particular case of the
projection $p\,:\,M\to pt$ gives the Lefschetz-type formula :
$$\forall s\in{\pi_{_M}}_*\Theta (p),\quad p_*(s)=
p^g_*(s\cdot e_{_{T_{M^g}M}}^{^{-1}})\leqno{\bf (2.9.3)}$$
where $s\cdot e_{_{T_{M^g}M}}^{^{-1}}$ is a meromorphic section of
${\pi_{_{M^g}}}_*\Theta(p^g)$ without poles on ${\cal X}_g\i{\cal X}_G$.

\vskip.3cm

\noindent{\bf (2.10) Localization of correspondences.}
Suppose that $G$ is a compact Abelian Lie group.
Let $M_1,M_2$, two smooth complex algebraic $G$-varieties and
$W$ a correspondence between $M_1$ and $M_2$, i.e. a $G$-stable
(possibly singular) subvariety $W\i M_1\times M_2$ such that
the second projection, $q_{_2}$, is proper.
{}From formula (2.4.1) we get
$$\Theta(q_{_2})=\Theta_W(M_{_1}\times M_{_2})\otimes
\Theta(q^*_{_1}TM_{_1})^{^{-1}}.$$
For any $g\in G$, denote by $i$ the inclusion $W^g\hookrightarrow W$.
The contravariant functoriallity of ${\rm Ell}_G$ with
respect to the inclusion of pairs
$(M_{_1}^g\times M^g_{_2},\, M^g_{_1}\times M^g_{_2}\setminus W^g)
\i (M_{_1}\times M_{_2},\, M_{_1}\times M_{_2}\setminus W)$
gives a morphism of sheaves over $W_{_{{\cal X}_G}}$
$$i^*\,:\,\Theta_W(M_{_1}\times M_{_2})\to
{i_{_{{\cal X}_G}}}_*\Theta_{W^g}(M^g_{_1}\times M^g_{_2}).$$
It induces a morphism of sheaves
$$\Theta(q_{_2})\to
{i_{_{{\cal X}_G}}}_*\Theta_{W^g}(M_{_1}^g\times M_{_2}^g)
\otimes\Theta(q^*_{_1}TM_1)^{^{-1}}=
{i_{_{{\cal X}_G}}}_*\bigl(\Theta(q_{_2}^g)\otimes
\Theta({q^g_{_1}}^*T_{M_1^g}M_1)^{^{-1}}\bigr).$$
The Euler class $e_{_{T_{M^g_1}M_1}}$
is a regular section of $\Theta(T_{M^g_1}M_1)^{^{-1}}$.
Thus, multiplying  by $e_{_{T_{M^g_1}M_1}}^{^{-1}}$ gives a rational
morphism of sheaves
$r_g\,:\,{\pi_{_{W}}}_*\Theta(q_{_2})\to{\pi_{_{W^g}}}_*\Theta(q_{_2}^g)$.

\proclaim (2.10.1) Proposition. For any $g\in G$, the map
$r_g\,:\,{\pi_{_{W}}}_*\Theta(q_{_2})\to{\pi_{_{W^g}}}_*\Theta(q_{_2}^g)$
is regular and invertible over ${\cal X}_g$.
Moreover $r_g$ commutes with the product of correspondences, as defined
in (2.7.2).

\noindent{\sl Proof.} The first part of the proposition follows from
(2.9). With the notations of (2.7.2), the second part
is reduced to the commutativity of the following diagram
$$\matrix{
{\pi_{_{W_{12}}}}_*\Theta(q_{_2})\otimes {\pi_{_{W_{23}}}}_*\Theta(q_{_3})&
\buildrel r_g\otimes r_g\over \longrightarrow&
{\pi_{_{W_{12}^g}}}_*\Theta(q^g_{_2})\otimes{\pi_{_{W^g_{23}}}}_*
\Theta(q^g_{_3})\cr
\downarrow&&\downarrow\cr
{\pi_{_{W_{13}}}}_* \Theta(\overline q_{_3})&
\buildrel r_g\over\longrightarrow&
{\pi_{_{W^g_{13}}}}_*\Theta(\overline q^g_{_3})
}, \leqno {\bf (2.10.2)}$$
whose proof is done exactly as in equivariant K-theory
(see [CG; Theorem 4.10.12]).

\vskip.2cm

\noindent{\bf (2.10.3)} Similarly, if $N_1=T^*M_1$, $N_2=T^*M_2$, and
$Z\i N_1\times N_2$ is a
Lagrangian correspondence, the Proposition (2.10.1) has the following analogue,
with the notations of (2.8) :
contravariant functoriality with respect to the inclusion, $i$, of
fixed points subvarieties gives a map
$$i^*\,:\,\Theta_Z(N_{_1}\times N_{_2})\to
{i_{_{{\cal X}_G}}}_*\Theta_{Z^g}(N^g_{_1}\times N^g_{_2}).$$
By composing it with the product by
$e_{_{T_{M^g_1}M_1}}^{^{-1}}\otimes e_{_{T_{M^g_2}M_2}}^{^{-1}}$,
we get a rational morphism of sheaves
$\rho_g\,:\,{\pi_{_{Z}}}_*\Xi_{_Z}\to{\pi_{_{Z^g}}}_*\Xi_{_{Z^g}}$,
which is regular and invertible over ${\cal X}_g$ and commutes with the
product of correspondences.

\vskip 1cm

\centerline {\bf  3. Geometric construction of current algebras.}

\vskip 1cm

\noindent {\bf (3.1)} Let $X$ be a smooth irreducible complex algebraic
variety. For any $d\geq 0$ we have a natural projection
$$\pi : X^d \rightarrow X^{(d)} = X^d/S_d,$$
where $X^{(d)}$ is the $d$-fold symmetric product of $X$. For any
coherent ${\cal O}_X$-sheaf $\cal F$  we construct coherent sheaves
${\cal F}^{(\oplus d)}$ and ${\cal F}^{(\otimes d)}$ on $X^{(d)}$
as follows :
$${\cal F}^{(\oplus d)} = \biggl(\pi_* \biggl(\bigoplus_{i=1}^d {\rm pr}_i^*
{\cal F}\biggl)\biggl)^{S_d}, \quad
{\cal F}^{(\otimes d)} = \biggl(\pi_* \biggl(\bigotimes_{i=1}^d {\rm pr}_i^*
{\cal F}\biggl)\biggl)^{S_d},$$
where ${\rm pr}_i: X^d\rightarrow X$ is the $i$th projection.
Thus the geometric fiber
of ${\cal F}^{(\oplus d)}$ (resp. ${\cal F}^{(\otimes d)}$) at a point
$I = \{x_1, ..., x_d\} \in X^{(d)}$ is given by $\bigoplus_{x\in I} {\cal F}_x$
(resp. $\bigotimes_{x\in I} {\cal F}_x$).

\vskip .3cm

\noindent {\bf (3.2)} Let $\cal G$ be a locally free ${\cal O}_X$-sheaf
of Lie algebras (with bracket not necessarily ${\cal O}_X$-linear). Then
${\cal G}^{(\oplus d)}$ has an obvious structure of a sheaf of Lie algebras
on $X^{(d)}$. Given a locally free ${\cal O}_X$-sheaf $\cal M$ of $\cal G$-
modules, we find, in the same way, that ${\cal M}^{(\otimes d)}$ is a sheaf
of ${\cal G}^{(\oplus d)}$-modules so that there is a natural Lie
algebra morphism
$${\cal G}^{(\oplus d)} \to {\cal E}nd {\cal M}^{(\otimes d)}\leqno {\bf
(3.2.1)}$$

For instance, if $X$ is a point, then $\cal G$ is a finite-dimensional
Lie algebra, $\cal M$ is a finite-dimensional $\cal G$-module, ${\cal G}^{
(\oplus d)} = {\cal G}$, and ${\cal M}^{(\otimes d)} = S^d{\cal M}$.

Observe that for any open set $U\i X$ we have the ``diagonal"
homomorphism of Lie algebras $\Gamma(U, {\cal G}) \rightarrow
\Gamma(U^{(d)}, {\cal G}^{(\oplus d)})$. Hence the space $\Gamma(U^{(d)},
{\cal M}^{(\otimes d)})$ acquires a natural structure of $\Gamma(U, {\cal G})$-
module.

We will be mainly concerned with the case ${\cal G} = {\cal E}nd ({\cal M})$,
where $\cal M$ is (the sheaf of sections of) a vector bundle on $X$.
For any open $U\in X$ the above construction makes
$\Gamma(U^{(d)},{\cal M}^{(\oplus d)})$ into a representation of
the ``current algebra" $\Gamma(U, {\cal E}nd ({\cal M}))$.

\vskip .3cm

\noindent {\bf (3.3)} We recall the general formalism of sheaf-theoretic
correspondences. Let $s: V\rightarrow S$ be a finite morphism of
normal algebraic varieties. Form the Cartesian product $M = V\times_S V$
and let $\tilde M$ be the normalization of $M$. We have the natural
commutative diagram of projections
$$\matrix{&\tilde M&\buildrel q_2\over\longrightarrow& V&\cr
q_1&\big\downarrow&\buildrel r\over\searrow&\big\downarrow&s\cr
&V&\buildrel s\over\longrightarrow&S&}.$$
We view $\tilde M$ as a correspondence from $V$ to $V$.
Let $\cal V$ be a vector bundle on $V$. We define the following sheaves on $S$
:
$${\cal W}=s_*{\cal V},\quad\tilde{\cal G}=
r_*{\cal H}om(q_2^*{\cal V},q_1^*{\cal V}).$$
The composition of ${\cal H}om$
makes $\tilde{\cal G}$ into a sheaf of associative algebras on $S$,
and $\cal W$ is a sheaf of $\tilde{\cal G}$-modules.

\vskip.2cm

\noindent{\bf (3.3.1) Remark.} The structure of $\tilde{\cal G}$-module on
$\cal W$ can be described in terms similar to those in (2.7.1). Indeed, fix
$s$ and $g$ two local sections of $\tilde{\cal G}$ and ${\cal W}$. They can
be seen as sections of ${\cal H}om(q_2^*{\cal V},q_1^*{\cal V})$ and
${\cal V}$. Lift $g$ to a section of $q_1^*{\cal V}$ and apply $s$ to it,
we get a section $s(g)$ of $q^*_2{\cal V}$. Now, since $q_2$ is flat and finite
and ${\cal V}$ is locally trivial we have a {\sl trace} map
${q_2}_*q_2^*{\cal V}\rightarrow{\cal V}$; just take the image of $s(g)$ by
this trace. This way we obtain a Lie algebra map
$$\tilde{\cal G}\rightarrow {\cal E}nd({\cal W})
\leqno {\bf (3.3.2)}$$
which is clearly injective by construction.

\vskip .3cm

\noindent {\bf (3.4)} Let now $q: Y\rightarrow X$ be an unramified covering
of a smooth irreducible variety $X$ ( $Y$ may be disconnected),
${\cal L}$ a line bundle on $Y$ and ${\cal M} = q_*{\cal L}$. In the setup of
the previous section, we put $S = X^{(d)}$, $\,V = Y^{(d)}$,
$\,s = q^{(d)}: Y^{(d)}\rightarrow X^{(d)}$ the map induced by $q$.
Put also ${\cal V} = {\cal L}^{\otimes d}$, ${\cal G} = {\cal E}nd({\cal M})$,
so that ${\cal W} = {\cal M}^{(\otimes d)}$ and
$$\tilde{\cal G} = r_* {\cal H}om (q_2^*{\cal L}^{\otimes d},
q_1^*{\cal L}^{\otimes d})\hookrightarrow {\cal E}nd({\cal M}^{\otimes d}).$$
The normalization  $\tilde M$ of $M = Y^{(d)}\times_{X^{(d)}} Y^{(d)}$
is in this case nothing  but $(Y\times_X Y)^{(d)}$, the symmetric power
of the smooth variety $Y\times_X Y$.
The ${\cal G}^{(\oplus d)}$ - action on
${\cal M}^{(\otimes d)}$ defined in (3.2), gives a map
${\cal G}^{(\oplus d)} \rightarrow {\cal E}nd ({\cal M}^{(\otimes d)})$.

\proclaim (3.5) Theorem. There is a natural Lie algebra
homomorphism
$\tau: {\cal G}^{(\oplus d)} \rightarrow \tilde{\cal G}$ of
sheaves
on $X^{(d)}$ making the
following diagram commute:
$$\matrix{{\cal G}^{(\oplus d)}\buildrel{(3.2.1)}\over\longrightarrow &
{\cal E}nd\,{\cal M}^{(\otimes d)}\cr
\big\downarrow&\nearrow\,(3.3.2)&\cr
\tilde{\cal G}&&}.$$

\vskip 2mm
Let ${\cal U}(L)$ denote the universal
enveloping algebra of  Lie algebra $L$. Applying the ``diagonal" map to $\tau$,
we get the first part of the following proposition :

\proclaim (3.6) Proposition. (a) For any open  $U\i X$, we have a
Lie algebra  homomorphism
 $$\tau_U: \Gamma(U, {\cal G}) \rightarrow \Gamma(U^{(d)},
\tilde{\cal G}).$$
(b) The following homomorphism of associative algebras
$$ {\cal U}(\Gamma(U, {\cal G})) \rightarrow \Gamma(U^{(d)},
\tilde{\cal G})$$
induced by $\tau_U$, is surjective.

Proof of Theorem 3.5 will be given in (3.11), the proof of the second
part of Proposition 3.6 in (3.14).

\vskip .3cm

\noindent {\bf (3.7) Partitions and matrices.}
We need to introduce some combinatorial notation, to
 be used throughout the paper.

Fix an integer $d\geq 1$.
Let ${\bf V}\subset{\bf N}^n$ be the set of
all partitions ${\bf v} = (v_1,\ldots,v_n)$
of $d$ into $n$ summands. For ${\bf v}\in {\bf V}$ let $[{\bf v}]_i
\i \{1, ..., d\}$
be the  $i$-th segment of the partition {\bf v}, i.e.,
the segment $[v_1 + ... + v_{i-1}+1, \, v_1 + ... + v_i]$. Let
$S_{\bf v} = S_{_{[{\bf v}]_1}}\times ... \times S_{_{[{\bf v}]_n}}$ be
the Young subgroup
in the symmetric group $S_d$, corresponding to {\bf v}, i.e., the
subgroup of permutations preserving each segment $[{\bf v}]_i$.

Let {\bf M} be the set of  the $(n\times n)$-matrices
 $A=(a_{ij})\in{\bf Z}_+^{n^2}$ with non-negative integral
entries such that $\sum_{i,j} a_{ij} = d$.
For  ${\bf v}^1, {\bf v}^2 \in {\bf V}$ we denote
$${\bf M}({\bf v}^1, {\bf v}^2) = \biggl\{ A\in {\bf M}:
\sum_j a_{ij} = v_i^1, \,\, \sum_i a_{ij} = v_j^2\biggl\},$$
so that ${\bf M}  = \coprod_{{\bf v}^1, {\bf v}^2 \in {\bf V}}
 {\bf M}({\bf v}^1, {\bf v}^2)$. For any given ${\bf v}^1, {\bf v}^2$
we have a natural map
$${\bf m} = {\bf m}_{{\bf v}^1, {\bf v}^2}: S_d \longrightarrow
{\bf M}({\bf v}^1, {\bf v}^2)$$
which identifies ${\bf M}({\bf v}^1, {\bf v}^2)$ with the double coset
space $S_{{\bf v}^1}\backslash S_d/S_{{\bf v}^2}$. Namely, ${\bf m}$
associates to a permutation $\sigma\in S_d$ the matrix ${\bf m}^\sigma$
with
$${\bf m}^\sigma_{ij} = {\rm Card} \, \bigl\{ \alpha\in [{\bf v}^1]_i:
\, \sigma(\alpha)\in [{\bf v}^2]_j\bigl\}.$$
Let us introduce a partial order $\preceq$ on ${\bf M}({\bf v}^1, {\bf v}^2)$
as follows. For $A = (a_{ij})$
and $B = (b_{ij})$ in ${\bf M}({\bf v}^1, {\bf v}^2)$ we say that
$A \preceq B$, if for any $1\leq i<j\leq n$ we have
$$\sum_{_{r\leq i\,;\,s\geq j}}a_{rs}\leq\sum_{_{r\leq i\,;\,s\geq j}}b_{rs}$$
and for any $1\leq j<i\leq n$ we have
$$\sum_{_{r\geq i\,;\,s\leq j}}a_{rs}\leq\sum_{_{r\geq i\,;\,s\leq j}}b_{rs}.$$
The map
${\bf m}: S_d\rightarrow {\bf M}({\bf v}^1, {\bf v}^2)$ is monotone with
respect to the Bruhat order (denoted $\leq$)
on the symmetric group; in other words, $A\leq B\Rightarrow A\preceq B$
for any $A,B\in{\bf M}({\bf v}^1,{\bf v}^2)$ (see [BLM , Lemma 3.6]).
The order $\preceq$ on ${\bf M} = \coprod {\bf M}({\bf v}^1, {\bf v}^2)$
is the disjoint union of the orders on the components
${\bf M}({\bf v}^1, {\bf v}^2)$, i.e., the elements of different
components are set to be incomparable.

\vskip .3cm

\noindent {\bf (3.8) The case ${\cal M} = {\cal O}_X^n$.}
We are now going to consider a special case of the construction of (3.4-6).
Let $X$ be a smooth irreducible algebraic variety.
Recall that
$X^{(m)}$ stands for the $m$-th symmetric power. For ${\bf v}\in{\bf V}$ set
$X^{({\bf v})} = \prod_i X^{(v_i)}$. Similarly, given
$A\in{\bf M}({\bf v}^1, {\bf v}^2)$ we  set
$\,X^{(A)} = \prod_{i,j} X^{(a_{ij})}\,$. so that one has
natural projections
$$X^{({\bf v}^1)} \buildrel{pr_1}\over\leftarrow
X^{(A)}\buildrel{pr_2}\over\rightarrow X^{({\bf v}^2)}\leqno {\bf
(3.8.0)}$$
We will also frequently use the following natural maps
$$\pi_{_{\bf v}}: X^{({\bf v})}\rightarrow X^{(d)}, \quad
\pi_{_A}: X^{(A)}\rightarrow X^{(d)}, \quad
p_{ij}: X^{(A)}\rightarrow X^{(a_{ij})}$$

\vskip .2cm

Put
$[n] = \{1,2, ..., n\}$. Let $Y = X\times [n]$ and $q: Y\rightarrow X$
be the natural projection. Put ${\cal L} = {\cal O}_Y$, so that
${\cal M} = {\cal O}_X^n$. Thus ${\cal G} = gl_n({\cal O}_X)$
is the sheaf of maps $X\rightarrow gl_n({\bf C})$ and $\tilde{\cal G}=r_*
{\cal O}_{\tilde M}$.

Observe that
$Y^{(d)} = \coprod_{{\bf v}\in {\bf V}} X^{(\bf v)}$. Hence, we get the
canonical direct sum decomposition
$${\cal M}^{(\otimes d)} = q^{(d)}_* {\cal O}_{Y^{(d)}} =
\bigoplus_{\bf v} \pi_{{\bf v}*} {\cal O}_{X^{({\bf v})}}.\leqno {\bf
(3.8.1)}$$

\noindent Following (3.3), we put
$$M = Y^{(d)}\times_{X^{(d)}}Y^{(d)} = \coprod_{{\bf v}^1, {\bf v}^2}
M_{{\bf v}^1, {\bf v}^2}, \quad
M_{{\bf v}^1, {\bf v}^2} = X^{({\bf v}^1)}\times_{X^{( d)}}X^{({\bf v}^2)}.$$
The variety $M_{{\bf v}^1, {\bf v}^2}$ is connected but reducible.
The irreducible components  of this variety
are parametrized by $n\times n$ matrices
$A = \|a_{ij}\| \in {\bf M}({\bf v}^1, {\bf v}^2)$. For such a matrix $A$
the irreducible component $M_A$ is the closure of the set
$$\biggl\{ \bigl( (I_1, ..., I_n), (J_1, ..., J_n)\bigl) \in X^{({\bf v}^1)}
\times X^{({\bf v}^2)}\biggl| \,\,\bigcup_\nu I_\nu = \bigcup_\nu J_\nu, \,\,
|I_\mu \cap J_\nu| = a_{\mu\nu}\biggl\}.\leqno {\bf (3.8.2)}$$
The normalization of the variety $M_A$ is nothing but $X^{(A)}$, so
the normalization of $M$ is
$$\tilde M = (Y\times_X Y)^{(d)} = (X\times [n]^2)^{(d)} =
\coprod_{A\in {\bf M}} X^{(A)}.$$

\vskip .2cm

For $1\leq i,j \leq n$ let $E_{ij} \in {\rm Mat}_n ({\bf Z})$ be the standard
matrix unit ($1$ at the spot $(i,j)$ and $0$ elsewhere). For every partition
${\bf v} = (v_1, ..., v_n)$ of $d$ and $i\neq j$ consider the following
integral $n\times n$-matrix
$$A_{ij}({\bf v}) = {\rm diag}(v_1, ..., v_{j-1}, v_j - 1, v_{j+1}, ..., v_n)
+ E_{ij}.$$
Let $\tilde M^{ij} = \coprod_{{\bf v}\in {\bf V}}
X^{(A_{ij}({\bf v}))}$. The projections $p_{ij}$ of the components give a map
$\tilde p_{ij}:
\tilde M^{ij} \rightarrow X$. For an open set $U\i X$ and
$f\in {\cal O}(U)$ we denote by $p_{ij}^*f$ the function $f\circ\tilde p_{ij}$
on $\tilde p_{ij}^{-1}(U)$ regarded as a function on the whole
$\coprod_{_{A\in{\bf M}}}U^{(A)}\i\tilde M$
(so this function is zero on components not in $\tilde p_{ij}^{-1}(U)$).
Similarly we denote by
$p_{ii}^*(f)$ the  function on $\coprod_{_{A\in{\bf M}}}U^{(A)}\i\tilde M$
which on components not of the form $U^{({\rm diag}({\bf v}))}$
is equal to 0 and on the component $U^{({\rm diag}({\bf v}))} =
U^{({\bf v})}$ is
$$p_{ii}^*(f)(I_1, ..., I_n) = \sum_{x\in I_i} f(x).$$

\noindent
We denote also by $E_{ij}(f)$ the section $E_{ij}\otimes f$ of ${\cal G} =
gl_n({\cal O}_X)$.

\proclaim (3.9) Theorem. The assignment
$E_{ij}(f) \mapsto p_{ij}^* f$
yields a Lie algebra homomorphism $\tau_U: \Gamma(U, {\cal G})
\rightarrow \Gamma(U^{(d)}, \tilde {\cal G})$
making the following diagram commute :
$$\matrix{
& \Gamma(U, {\cal G})\buildrel{(3.3.1)}\over\longrightarrow
&\Gamma(U^{(d)},
{\cal E}nd\,{\cal M}^{(\otimes d)})\cr
\tau_U&\big\downarrow&\nearrow\,(3.3.2)&\cr
&\Gamma(U^{(d)}, \tilde{\cal G})&&}.$$

\noindent {\sl Proof.} Since the map $\Gamma(U^{(d)},\tilde{\cal G})
\rightarrow \Gamma(U^{(d)}, {\cal E}nd
{\cal M}^{(\otimes d)})$ is injective by (3.3), it is enough to show
that the image of each $E_{ij}(f)$ under $\tau_U$ comes from an element
of $\Gamma(U^{(d)}, \tilde{\cal G})$.

Consider the standard $gl_n({\bf C})$-action on
the space ${\bf C}^n$ with standard basis $e_1, ..., e_n$. The
element $E_{ij}$ takes $e_i\mapsto e_j$ and other base vectors to 0.
Let $I$ be any finite set, $|I|=d$, and let ${\cal O}(I)$ be the ring of
functions $I\rightarrow {\bf C}$. The Lie algebra $gl_n(I) =
gl_n({\bf C})\otimes {\cal O}(I)$ acts on the tensor product of $I$ copies
of ${\bf C}^n$, which we denote by $({\bf C}^n)^{\otimes I}$. The natural basis
in $({\bf C}^n)^{\otimes I}$ is labelled by partitions ${\bf I} =
(I_1, ..., I_n)$, $I_\nu \i I,$ $I = \coprod I_\nu$. More precisely, the
basis vector corresponding to $I$ is $e_{\bf I} = \bigotimes_{\nu\in I}
e_{i(\nu)}$ where $i(\nu)$ is such that $\nu\in I_{i(\nu)}$.
For $f\in {\cal O}(I)$ and distinct $i,j\in [n]$ the standard action of
$E_{ij}(f) = E_{ij} \otimes f$ on  $({\bf C}^n)^{\otimes I}$ is immediately
seen to be given by the formula
$$e_{\bf I} \mapsto \sum_{\nu\in I_j} f(\nu) e_{(I_1, ..., I_i+\{\nu\}, ...,
I_j - \{\nu\}, ..., I_n)}.\leqno {\bf (3.10)}$$
The assignment $E_{ij}(f)\mapsto p_{ij}^*f$ in Theorem 3.8 is
nothing but the result of the simultaneous application of (3.10) to all subsets
$I\i U$, $|I|=d$. This completes the proof.

\vskip .3cm

\noindent {\bf (3.11) Proof of Theorem 3.5.} We first consider the case when
the covering $Y\rightarrow X$ and the bundle $\cal L$ on $Y$ are trivial, i.e.,
$Y=X\times [n]$, ${\cal L}={\cal O}_Y$, so that we are in the
situation of (3.7).
Let $U\i E$ be an affine open set. Then the space
$$\Gamma(U^{(d)}, {\cal G}^{(\oplus d)}) = \biggl(\bigoplus_{i=1}^d
\Gamma(U^d, {\rm pr}_i^*{\cal G})\biggl)^{S_d}$$
is generated, as a module over $\Gamma(U^{(d)}, {\cal O})$, by the image of
$\Gamma(U, {\cal G})$ under the diagonal map, i.e., by elements of the form
$(g(x_1), ..., g(x_n))$, $g\in \Gamma(U, {\cal G})$. We claim that the
assignment $\tau_U: E_{ij}(f)\mapsto p_{ij}^* f$
of Theorem 3.8 extends, by ${\cal O}_{U^{(d)}}$-linearity, to a morphism
$\Gamma(U^{(d)}, {\cal G}^{(\oplus d)})\rightarrow \Gamma(U^{(d)}, \tilde
{\cal G})$. Indeed, since ${\cal G} = \bigoplus_{i,j} {\cal O}_X \cdot E_{ij}$
as a coherent sheaf, we find that ${\cal G}^{(\oplus d)} =
\bigoplus_{i,j} {\cal O}_X^{\oplus d}\cdot E_{ij}$. Thus
it is enough to consider every ${\cal O}_X\cdot E_{ij}$ separately
and to show that the following statement holds :

\proclaim (3.12) Lemma. Let $\phi_\alpha(x_1, ..., x_d) \in {\cal O}(U^d)$ be
symmetric functions and $f_\alpha \in\Gamma(U, {\cal O})$ be such that
$$\sum_\alpha \phi_\alpha(x_1, ..., x_d) (f_\alpha(x_1), ..., f_\alpha(x_d)) =
0 \quad {\rm in} \quad {\cal O}(U^d) \oplus ... \oplus {\cal O}(U^d).$$
Then the section $\sum \phi_\alpha(x_1, ..., x_d) p_{ij}^* f_\alpha \in$
$\Gamma(U^{(d)}, \tilde{\cal G})$ is equal to 0 for any $i,j$.

\noindent We leave the verification of this lemma to the reader.
Theorem 3.5 follows in
the case $Y=X\times [n]$, ${\cal L}={\cal O}_Y$.

To treat the general case, let $n$ be the degree (number of sheets) of the
map $Y\rightarrow X$. Choose an etale  covering $\{\psi_\alpha:
U_\alpha \rightarrow X\}$ (so we have  the maps
$p_\alpha: \psi_\alpha^* Y   =  Y\times_X U_\alpha \rightarrow Y$)
and choose trivializations
$k_\alpha: u_\alpha\times [n]\rightarrow \psi_\alpha^*Y$ and $k_\alpha^*
p_\alpha^*{\cal L}\rightarrow {\cal O}_{U_\alpha\times [n]}$.
Then apply the assignments $\tau_{U_\alpha}$ of Theorem 3.8 and check that they
satisfy the
descent condition to get a map defined over all the $X^{(d)}$. We leave the
easy details to the reader.

\vskip .3cm

\noindent {\bf (3.13) Algebra of convolution operators.}
We consider the situation of (3.8). Fix an open
set $U\i X$. Geiven a matrix $A \in {\bf M}$ we have the following
diagdam, see (3.8.0):
$$X^{({\bf v}^1)} \buildrel{pr_1}\over\leftarrow
X^{(A)}\buildrel{pr_2}\over\rightarrow X^{({\bf v}^2)}$$
To any function $f\in{\cal O}(U^{(A)})$ we associate the convolution
operator
$$\Delta(A,f): {\cal O}(U^{(\bf v^1)}) \to {\cal O}(U^{(\bf v^2)})
\quad,\quad
\Delta(A,f): \psi \mapsto (pr_2)_*((pr^*_1\psi)\cdot f)$$
Given a linear operator $F$ of the form
$F={\Delta}(A,f)+\sum_i{\Delta}(B^i,g^i)$, we write
$F={\Delta}(A,f)+\cdots$, if all $B^i < A$
with respect to the Bruhat order (3.7).

We start by describing, in abstract terms, the composition of operations
of type $\Delta(A, f)$. Let {\bf T} be the set of 3-dimensional arrays
$(t_{ijk})_{1\leq i,j,k\leq n}$ of non-negative integers,  such that
$\sum_{i,j,k} t_{ijk} = d$.
For any ${\bf v}^1, {\bf v}^2, {\bf v}^3 \in {\bf V}$ put
$${\bf T}({\bf v}^1, {\bf v}^2, {\bf v}^3) =
\biggl\{  T\in {\bf T}: \, \sum_{j,k} t_{ijk} = v^1_i,\,\,
\sum_{i,k} t_{ijk} = v^2_j, \,\, \sum_{i,j} t_{ijk} = v^3_k\biggl\}.$$
If $1\leq l\leq m \leq 3$ and $T\in {\bf T}({\bf v}^1, {\bf v}^2, {\bf v}^3)$,
let $T^{lm}$ be the matrix in ${\bf M}({\bf v}^l, {\bf v}^m)$
obtained by summing the entries of $T$ along the indexed not named.
For $T\in {\bf T}$ we set $X^{(T)} = \prod_{i,j,k} X^{(t_{ijk})}$ and write
$pr_{lm} : X^{(T)} \rightarrow X^{(T^{lm})}$ for the natural
projection.
Given two matrices $A\in{\bf M}({\bf v}^1,{\bf
v}^2)$, and
$B\in{\bf M}({\bf v}^2,{\bf v}^3)$, put
$${\bf T}(A,B):=\{T\in{\bf T}: \, T^{12}=A,
T^{23}= B\}\eqno {\bf (3.13.0)}$$
The definition of the convolution gives (3.3):

\proclaim (3.13.1) Lemma. Let $A\in{\bf M}({\bf v}^1,{\bf
v}^2)$, and
$B\in{\bf M}({\bf v}^2,{\bf v}^3)$, $f\in{\cal O}(U^{(A)})$,
$g\in {\cal O}(U^{(B)})$.
 Then, the composition of convolution operators on the LHS below
is given by
the following sum on the RHS:
$$\Delta(A,f)\cdot\Delta(B, g)=\sum_{\{T\in {\bf T}(A,B)\}}
 \Delta\bigl(T^{^{13}},\,pr_{_{13*}}
(pr_{_{12}}^*(f)\cdot pr_{_{23}}^*(g))\bigr).$$

\proclaim (3.13.2) Corollary. Suppose
that, in the situation of Proposition 3.13.1, $A = {\rm diag} (d^1_1,
..., d^1_n)$
is the  diagonal matrix so that ${\bf v}^1 = {\bf v}^2$.
Then we
 have
$$\Delta(A,f)\, \Delta(B,g) = \Delta(B, h),$$
where $h = g\cdot p_2^*f$ and $p_2: U^{(B)} \rightarrow U^{({\bf v}^2)} =
U^{(A)}$ is the natural projection.

\noindent {\sl Proof.} The only array $T$ which will contribute  to the sum
in (3.13.1), is given by
$$t_{ijk} = \cases{ 0, \quad i\neq j\cr b_{jk}, \quad i=j}$$

The following result gives way to proving results by induction
with respect to the Bruhat order.

{\proclaim (3.13.3) Proposition. In the situation of
Proposition 3.13.1, assume that
 $A={\rm diag} + a_{p,p+1}E_{p,p+1}$ where `diag' stands for  diagonal matrix.
Put $m = \max\{i: b_{p+1, i}\neq 0\}$, and suppose, in addition,
that $b_{p+1,m}\geq a_{p,p+1}$.  Then $$\Delta(A,f)
\cdot \Delta(B, g) = \Delta(C, h) + ..., $$
where $C = B + a_{p, p+1}(E_{pm} - E_{p+1, m})$ and
$h = pr_{13*}(pr_{12}^*(f) \cdot pr_{23}^*(g))$,
the projections $pr_{ij}$ corresponding to the array in (3.13.4).

Similarly, if $A - a_{p, p-1}E_{p, p-1}$ is diagonal,
 $m = \min\{ i: b_{p-1, i}\neq 0\}$ and $b_{p-1,m}\geq a_{p,p-1}$,
 then we have a similar relation with $C = B + a_{p, p-1}(E_{pm}-
E_{p-1, m})$
 and the array $T$ defined in an analogous way.}

To prove the proposition, we need the following result.

\proclaim (3.13.4) Lemma. Let
 $A={\rm diag} + a_{p,p+1}E_{p,p+1}$. For an arbitrary matrix
$B$ as above, the set ${\bf T}(A,B)$, see $(3.13.0)$, is in bijective
correspondence with
the following set of $n$-tuples
$${\bf S}(A,B)=\{s=(s_1,s_2,...,s_n)\,
|\, 0\leq s_k\leq b_{p+1,k}\quad ,\quad \sum_ks_k=a_{p,p+1}\}$$
To an $n$-tuple $s\i {\bf S}(A,B)$ on assigns the following array
$T(s)=(t_{ijk})$
$$
\cases{ t_{iik} = b_{ik} \quad {\sl for}\quad
i\neq p+1\cr t_{p+1, p+1, k} = b_{p+1, k} - s_k\cr t_{p, p+1, k} =
 s_k\cr t_{ijk} = 0 \quad {\sl otherwise,}}$$
Moreover, we have
$T(s)^{^{13}}=B+\sum_j s_j\cdot(E_{pj}-E_{p+1,j}).$

\noindent {\sl Proof of (3.13.4).} The arrays $T=(t_{ijk})$ which will
contribute to the sum in (3.13.1) are given by :
$$\cases{ \sum_kt_{ijk} = a_{ij}\cr \sum_it_{ijk} = b_{jk}.}$$
In particular,
$$\cases{ t_{ijk} = 0 \quad {\rm if}\quad (i,j)\neq (i,i),(p,p+1)\cr
t_{jjk} = b_{jk} \quad {\rm if}\quad j\neq p+1\cr
t_{p, p+1, k}+t_{p+1,p+1,k} = b_{p+1,k}.}$$
Thus $T$ will have a non-zero contribution if and only if :
$$\cases{ t_{ijk} = 0 \quad {\rm if}\quad (i,j)\neq (i,i),(p,p+1)\cr
t_{jjk} = b_{jk} \quad {\rm if}\quad j\neq p+1\cr
t_{p+1, p+1, k} = b_{p+1,k}-t_{p,p+1,k}\cr
t_{p, p+1, k}\leq b_{p+1,k}\cr
\sum_kt_{p, p+1, k}=a_{p,p+1}\cr
\sum_kt_{p+1, p+1, k}=a_{p+1,p+1}\cr
\sum_kb_{jk}=a_{jj} \quad {\rm if}\quad j\neq p+1.}$$
The last two equations follow directly from $\sum_jb_{ij}=d^2_i=\sum_ja_{ji}$.
The conditions above mean that there exists an $n$-tuple $s\i
{\bf S}(A,B)$ such that
$$\cases{ t_{iik} = b_{ik} \quad {\rm for}\quad
i\neq p+1\cr t_{p+1, p+1, k} = b_{p+1, k} - s_k\cr t_{p, p+1, k} =
t_k\cr t_{ijk} = 0 \quad {\rm otherwise.}}$$
The lemma follows.

\vskip.2cm

\noindent {\sl Proof of (3.13.3).} For $s\in{\bf S}(A,B)$ we put
$C(s):= T(s)^{^{13}}$. Observe that, if
$s,s'\in{\bf S}(A,B)$ then we have
$$C(s')\preceq C(s)\Leftrightarrow\forall j\quad,\quad\cases
{\sum_{r\geq j}s'_r\leq\sum_{r\geq j}s_r,\cr \sum_{r\leq j}s'_r\geq\sum_{r\leq
j}s_r.}$$
Thus, since $\sum_ks'_k=\sum_ks_k(=a_{p,p+1})$,
$$C(s')\preceq C(s)\Leftrightarrow
\sum_{r\geq j}s'_r\leq\sum_{r\geq j}s_r\quad,\quad\forall j.$$
Observe, on the other hand that
$b_{p+1,k}=0$ if $k>m$. Hence $s_k=0$ if $k>m$, for any
$s\in{\bf S}(A,B)$. Since $b_{p+1,m}\geq a_{p,p+1}$, the set
$\{T(s)\,|\,s\in{\bf S}(A,B)\}$ has a greatest
element with respect to $\preceq$; it corresponds to the $n$-tuple
$s=a_{p,p+1}\cdot(\delta_{1m},\delta_{2m},...,\delta_{nm})\in{\bf
T}(A,B)$.
The corresponding array looks like
$$T_{max}= \,\cases{ t_{ijk} = 0 \quad {\rm if}\quad (i,j)\neq (i,i),
(p,p+1)\cr
t_{iik} = b_{ik}-\delta_{_{(i,k)\,(p+1,m)}}a_{p,p+1}\cr
t_{p,p+1,k} = \delta_{km}a_{p,p+1}.}$$
Observe that the matix $C = T_{max}^{^{13}}$ arising from this array
is precisely the matrix $C=B+a_{p,p+1}(E_{pm}-E_{p+1,m})$.
In particular, for any other matrix $C'\neq C$ arising on the RHS
of the formula
$$\Delta(A,f)
\cdot \Delta(B, g) = \sum_{C',h} \Delta(C', h) + ..., $$
we have $C'\preceq C$. Thus, proposition
(3.13.3) is proved.

\vskip .2cm

\noindent {\bf (3.14) Proof of Proposition 3.6 (b).}
We will prove by induction on
$l(C)=\sum_{i\neq j}\pmatrix{|i-j|+1\cr 2} c_{ij}$
that, for any $C\in{\bf M}({\bf v}^1,{\bf v}^2)$
and $h\in{\cal O}(U^{(C)})$, the convolution
operator $\Delta(C,h)$ belongs to ${\rm Im}\,{\cal U}(\tau_U)$. If $l(C)=0,1$
it is evident. Suppose $l(C)>1$. Then $C$ is not diagonal. Suppose for
instance that $C$ is not a lower-triangular matrix. Let $(p,q)$
be the greatest element in
$\{(i,j)\,|\,1\leq i<j\leq n\quad,\quad c_{ij}\neq 0\}$
with respect to the right-lexicographic order (by right-lexicographic
order we mean $(i,j)>(p,q)\Leftrightarrow j>q\quad{\rm or}\quad
(j=q\quad {\rm and}\quad i>p)$). Put
$$B=C+c_{pq}(E_{p+1,q}-E_{pq})\quad{\rm and}\quad
A={\rm diag}(v^1_1,v^1_2,...,v^1_n)+c_{pq}(E_{p,p+1}-E_{p,p}).$$
Then $l(B)<l(C)$, $q={\rm max}\{i\,|\,b_{p+1,i}\neq 0\}$,
$$b_{p+1,q}=c_{pq}+c_{p+1,q}\geq c_{pq}=a_{p,p+1},$$
and, for any
$f\in{\cal O}(U^{(A)})$, $g\in{\cal O}(U^{(B)})$, (3.13.3) gives :
$$\Delta(A,f)\cdot\Delta(B,g)=\Delta(C,h)+...$$
where $h=pr_{_{13*}}(pr^*_{_{12}}(f)\cdot pr^*_{_{23}}(g))$, the projections
$pr_{_{ij}}$ beeing computed with respect to the array $T$ given in (3.13.6).
In particular if $f=p^*_{_{p,p+1}}(f_{_{p,p+1}})$, $f_{_{p,p+1}}\in{\cal O}
(U^{(a_{p,p+1})})$, one gets :
$$h=p^*_{_{pq}}(f_{_{p,p+1}})\cdot r^*(g)\in{\cal O}(U^{(C)}),$$
where $r\,:\,U^{(C)}\rightarrow U^{(B)}$ is given by :
$$\cases{
r_{|U^{(c_{ij})}}=id\,:\,U^{(c_{ij})}\rightarrow U^{(b_{ij})} &if $(i,j)\neq
(p,q),(p+1,q)$,\cr
r_{|U^{(c_{pq})}\times U^{(c_{p+1,q})}}=\oplus\,:\,U^{(c_{pq})}\times
U^{(c_{p+1,q})}\rightarrow U^{(b_{p+1,q})}.}$$
It's then easy to conclude from the induction hypothesis since that
$\Delta(A,f)\in{\rm Im}\,{\cal U}(\tau_U)$.

\vskip 2cm

\centerline {\bf 4. Elliptic cohomology construction of current algebras.}

\vskip 1cm

\noindent {\bf (4.1) The Steinberg variety.} Fix $d,n \geq 0$. Let $F$
be the variety of $n$-step flags in ${\bf C}^d$, i.e., filtrations
$D = \{0=D_0\subseteq D_1\subseteq D_2 \subseteq ...
\subseteq D_n = {\bf C}^d\}$. The connected components of $F$
 are parametrized by partitions ${\bf v} = (v_1, ..., v_n)$ of $d$;
the component $F_{\bf v}$ corresponding to {\bf v} consisting of
flags $D$ such that $\dim D_i/D_{i-1} = v_i$. As in section
(3.7), we denote the set of such partitions {\bf V}; we also conserve the
other notations of that section.

Consider the group $GL_d({\bf C})$ acting diagonally on $F\times F$.

\proclaim  (4.1.1) Proposition. (a) The orbits of $GL_d({\bf C})$ are
parametrized by matrices $A\in {\bf M}$; the orbit $O_A$ corresponding to
$A$ consists of pairs of flags $(D, D')$ such that
$\dim (D_i\cap D'_j) = a_{ij}$. \hfill\break
(b) The closure of $O_A$ contains $O_B$ if and only if $B\leq A$
with respect to the Bruhat order on {\bf M}.

Denote by
$$Z\quad = \quad \bigcup_{A\in M} \overline{ T^*_{O_A} (F\times F)}
\quad \i \quad T^*(F\times F) = T^*F\times T^*F$$
the union of the closures of the conormal bundles to $O_A$ (the Steinberg-
type variety). We regard $Z$ as a correspondence from $T^*F$ to $T^*F$.
It is known that $Z$ satisfies the following two properties :

\vskip .2cm

\noindent {\bf (4.1.2)} Both projections $Z\rightarrow T^*F$
are proper.

\vskip .2cm

\noindent {\bf (4.1.3)} The composition $Z\circ Z$ is equal to $Z$.

\vskip .3cm

\noindent {\bf (4.2) Elliptic cohomology of the Steinberg variety.}
Fix an elliptic curve $E\rightarrow S$ over a base scheme $S$.

Put $G = U(d)$ and let $T\i G$ be the subgroup of diagonal matrices.
Denote by ${\rm Ell}_G$, ${\rm Ell}_T$
the $G$-equivariant and $T$-equivariant elliptic cohomology
theories with values in coherent sheaves on ${\cal X}_G = E^{(d)}$
and ${\cal X}_T = E^{d}$. We conserve all the other notations of Section 1.

We have an obvious $G$-action on $Z$ and $T^*F$ and the second projection
$Z\rightarrow T^*F$ is $G$-equivariant (and proper). Thus
$Z$ is a Lagrangian correspondence between $F$ and $F$, see (...)
and so we have
the microlocal Thom sheaf $\Xi_{_Z}$ on $Z_{{\cal X}_G}$. Since $Z\circ Z=Z$,
it follows from Proposition (2.8.6) that the direct image ${\pi_{_Z}}_*
\Xi_{_Z}$ is a sheaf of associative algebras on ${\cal X}_G$,
and ${\pi_{_F}}_*\Theta(T^*F)^{-1}$ is a sheaf of its modules.
On the other hand, ${\cal X}_G = E^{(d)}$, and in (3.3) we
constructed a sheaf of algebras $\tilde{\cal G} = \bigoplus_{A\in {\bf M}}
{\pi_{_A}}_* {\cal O}_{E^{(A)}}$ on $E^{(d)}$.

\proclaim (4.3) Theorem. There is an isomorphism of sheaves of algebras
$\tilde{\cal G}\rightarrow {\pi_{_Z}}_*\Xi_{_Z}$.

\noindent The rest of this section is devoted to the proof of Theorem 4.3.

\vskip .3cm

\noindent {\bf (4.4)}
We start by considering closed orbits in $F\times F$. Fix ${\bf v}\in {\bf V}$
and $i,j\in \{1, ..., n\}$ such that $|i-j|\leq 1$.
For such $i,j$ let
$A_{ij}({\bf v})$ be the matrix introduced in (3.8).
The orbit $O_{A_{ij}({\bf v})}$ is the incidence correspondence in
$F_{\bf v}\times F_{{\bf v}'}$, where ${\bf v}' = (v_1,..., v_i+j-i, v_{i+1}
-j+i,..., v_n)$. We have thus the diagram of flag varieties
$$F_{\bf v} \longleftarrow O_{A_{ij}({\bf v})}\longrightarrow F_{{\bf v}'}$$
and the result of application of the functor ${\cal X}_G$ to this diagram
is just
$$E^{({\bf v})} \leftarrow E^{(A_{ij}({\bf v}))} \rightarrow
E^{({\bf v}')},$$
see (1.7.5). Denote for short
$$Z_{ij}({\bf v})= T^*_{O_{A_{ij}({\bf v})}}(F_{\bf v}\times F_{{\bf v}'}).$$
This is an irreducible component of the ``Steinberg-type" variety $Z$.
Note that $T^*F_{\bf v}$, $T^*F_{{\bf v}'}$ and $Z_{ij}({\bf v})$
are $G$-equivariantly homotopy equivalent to $F_{\bf v}$, $F_{{\bf v}'}$
and $O_{A_{ij}({\bf v})}$ respectively; thus we can identify, in
particular, $Z_{ij}({\bf v})_{{\cal X}_G}$ with $E^{(A_{ij}({\bf v}))}$.
The right map in the diagram
$$T^*F_{\bf v} \longleftarrow Z_{ij}({\bf v})\longrightarrow T^*F_{{\bf v}'}$$
is proper and
induces, therefore, the  microlocal Thom sheaf $\Xi_{Z_{ij}({\bf v})}$
on $E^{(A_{ij}({\bf v}))}$. Proposition 2.8.6 implies the first
part of the following assertion :

\proclaim (4.5) Proposition. (a) If $|i-j|\leq 1$, the sheaf
$\Xi_{Z_{ij}({\bf v})}$ is naturally isomorphic to
${\cal O}_{E^{(A_{ij}({\bf v}))}}$.\hfill\break
(b) The sheaf $\Theta(T^*F)^{-1}$ on $F_{{\cal X}_G} = \coprod_{{\bf v}\in
{\bf V}}E^{({\bf v})}$ consists of functions $g(I_1, ..., I_n)$,
$I_\nu\in E^{(v_\nu)}$ which are allowed poles of at most 1-st order
along the locus of $(I_1, ..., I_n)$ such that $I_\mu\cap I_\nu\neq
\emptyset$ for some $\mu\neq \nu$.\hfill\break
(c) If $|i-j|\leq 1$, the map
$${\pi_{Z_{ij}({\bf v})}}_* \Xi_{Z_{ij}({\bf v})}
\rightarrow {\cal H}om ({\pi_{_{F_{\bf v}}}}_* \Theta(T^*F_{\bf v})^{-1},
{\pi_{_{F_{\bf v'}}}}_*\Theta(T^*F_{\bf v'})^{-1})$$
from (2.7.1) associates to $f\in {\cal O}(U^{(A_{ij}({\bf v}))})$,
$U\i E$, the convolution operator $\Delta(A_{ij}({\bf v}),f)$ from
(3.13).

\noindent {\sl Proof.} (b) Denote by ${\cal D}_1\i ... \i {\cal D}_n$ the
tautological
flag of bundles on $F$. It is well known that $TF$ has a $G$-equivariant
filtration whose quotients are
${\cal H}om ({\cal D}_i/{\cal D}_{i-1}, {\cal D}_j/{\cal D}_{j-1})$,
$i<j$. Our statement follows from this and the definition of $\Theta$.
Part (c) follows from the localisation formula given in (2.9) (see [V]).

\vskip .3cm

\noindent {\bf (4.6)} The rest of the proof of (4.3) is done by induction,
similarly to the reasoning of (3.13), from which we freely
borrow the notation. We use, in particular, the Bruhat order $\leq$
and the binary operation $*$ on the set ${\bf M}$. For
$A\in {\bf M}$, an open $U\i E$ and $f\in {\cal O}(U^{(A)})$ we denote
$\Delta(A,f)$ the section of $\tilde{\cal G}$ over $U^{(d)}$ corresponding
to $f$. We denote by
$\tilde{\cal G}_{\leq A}=\bigoplus_{_{B\leq A}}{\pi_{_B}}_*{\cal O}_{E^{(B)}}$
the subsheaf in
$\tilde{\cal G}$ spanned (over ${\cal O}_{E^{(d)}}$) by sections
of the form $\Delta(B,f)$ where $B\leq A$. Then
$\tilde{\cal G}_{\leq A}\cdot\tilde{\cal G}_{\leq B}
\i\tilde{\cal G}_{\leq A*B}$.
So we have a filtration of the sheaf of algebras $\tilde{\cal G}$
labelled by {\bf M}.

\vskip.2cm

On the other hand, let us denote the sheaf ${\pi_{_Z}}_*\Xi_{_Z}$ on
$E^{(d)}$ by $\widehat{\cal G}$. The Bruhat order, $\leq$, on ${\bf M}$ is,
in geometric terms, just the inclusion of orbit closures in $F\times F$
(i.e.  $\overline{O_A} = \bigcup_{B\leq A} O_B$). Let
$$Z_{\leq A} = \bigcup_{B\leq A}T^*_{O_B}(F\times F), $$
and $Z_{<A}$ be the similar union over $B<A$. The complement
$Z_A=Z_{\leq A}\setminus Z_{<A}=T^*_{O_A}(F\times F)$ is a smooth
locally-closed $G$-sub-variety of $T^*(F\times F)$ paved by complex affine
spaces. Let $i_{_{\leq A}}$ be the inclusion maps of $Z_{\leq A}$ in $Z$.
By induction, using the long exact sequence in ``homology" (2.2.1)
and the vanishing of ${\rm Ell}_G^{2k+1}({\bf C}^n)$,
we get that the Gysin map (see (2.8.8.2))
$${i_{_{\leq A}}}_*\,:\,{\pi_{_{Z\leq A}}}_*\Xi_{_{Z\leq A}}\rightarrow
{\pi_{_Z}}_*\Xi_{_Z}=\widehat{\cal G}$$
is an injective homomorphism of coherent sheaves on $E^{(d)}$.
Let $\widehat{\cal G}_{\leq A}$ be its image. Then
$\widehat{\cal G}_{\leq A}\cdot\widehat{\cal G}_{\leq B}\i
\widehat{\cal G}_{\leq A*B}$.
We have thus constructed a filtration of the sheaf of algebras
$\widehat{\cal G}$ labelled by ${\bf M}$.

\proclaim (4.6.1) Lemma. The sheaf $\,\widehat{\cal G}$ is locally free.
Moreover, the morphism
$\Psi\,:\,\widehat{\cal G}\to{\cal E}nd
\bigl(\pi_{_{F*}}\Theta(T^*F)^{^{-1}}\bigr)$
induced by the action of $\widehat{\cal G}$ on
$\pi_{_{F*}}\Theta(T^*F)^{^{-1}}$
is injective.

\noindent{\sl Proof.} The long exact sequence (2.2.1) gives short
exact sequences
$$ 0\to\widehat{\cal G}_{< A}\to\widehat{\cal G}_{\leq A}
\to{\pi_{_{Z_A}}}_*\Xi_{_{Z_A}}\to 0,$$
where $\Xi_{_{Z_A}}$ is defined as in (2.8.8.1).
Since
$\Xi_{_{Z_A}}\simeq{\cal O}_{{Z_A}_{{\cal X}_G}}\simeq{\cal O}_{E^{(A)}}$,
an induction on $A\in{\bf M}$ proves that $\widehat{\cal G}$ is locally free.
Set $E^{(d)}_{_{reg}}=E^{(d)}\setminus\bigcup_{H\i U(d)}{\cal X}_H$,
the union beeing taken over all the closed proper subgroups of $U(d)$.
For any regular element $g\in U(d)$, the morphism
of localisation of correspondences (see (2.10.3))
$$\rho_g\,:\,{\pi_{_{Z}}}_*\Xi_{_Z}\to{\pi_{_{Z^g}}}_*\Xi_{_{Z^g}}$$
is an isomorphism of sheaves of algebras over $E^{(d)}_{_{reg}}$.
But $Z^g$ is formed of a finite number of points of $Z$ labelled
by ${\bf M}$, and one can thus easily see that
the map $\Psi$ is an isomorphism over $E^{(d)}_{_{reg}}\i E^{(d)}$.
Since $\widehat{\cal G}$ is locally free, $\Psi$ is injective.

\vskip.2cm

We now proceed to construct a homomorphism of sheaves of algebras
$\Phi: \tilde{\cal G}\rightarrow \widehat{\cal G}$. Namely, from
(3.14) the sheaf of algebras $\tilde{\cal G}$ is generated by its subsheaves
$\tilde{\cal G}_{\leq A_{ij}({\bf v})}$ for ${\bf v}\in {\bf V}$
and $|i-j|\leq 1$. Such matrices
$A_{ij}({\bf v})$ are minimal with respect to $\leq$, since the
corresponding orbits are closed. Thus the arguments of (4.4)
give identifications $\Phi_{ij{\bf v}}:
\tilde{\cal G}_{\leq A_{ij}({\bf v})} \rightarrow
\widehat{\cal G}_{\leq A_{ij}({\bf v})}$.

\proclaim (4.6.2) Lemma. The isomorphisms  $\Phi_{ij{\bf v}}$
extend to a unique homomorphism $\Phi: \tilde{\cal G}\rightarrow
\widehat{\cal G}$ of sheaves of algebras such that
$\Phi(\tilde{\cal G}_{\leq A}) \i \widehat{\cal G}_{\leq A}$ for any $A$.

\noindent {\sl Proof.} The sheaf $\pi_{_{F*}}\Theta(T^*F)^{^{-1}}$
is locally free, since it is a
pushdown of a line bundle with respect to the finite flat map
$\pi_{_F}\,:\,F_{_{{\cal X}_G}}=
\coprod_{_{\bf v}}E^{({\bf v})}\rightarrow E^{(d)}$.
Therefore the sheaf of algebras
${\cal E}nd\bigl(\pi_{_{F*}}\Theta(T^*F)^{^{-1}}\bigl)$
on $E^{(d)}$ is locally free as well.

Let ${\bf R}$ be the algebra formed by operators
$\Delta(A,f)$ where $A\in {\bf M}$, and $f$
is a rational function on $E^{(A)}$.
Denote ${\cal R}$ the constant sheaf
on $E^{(d)}$ with stalk ${\bf R}$.
We claim that
${\cal E}nd\bigl(\pi_{_{F*}}\Theta(T^*F)^{^{-1}}\bigl)$
is naturally a subsheaf in ${\cal R}$.
Indeed, let $E^{(d)}_{_{gen}}\i E^{(d)}$
be the open set formed by unordered tuples of {\it distinct} points.
Then the map $\pi_{_F}$ is unramified over $E^{(d)}_{_{gen}}$.
Moreover, the line bundle
$\Theta(T^*F)^{^{-1}}$ is trivialized over
$\pi_{_F}^{^{-1}}(E^{(d)}_{_{gen}})$ by
Proposition 4.5(b). This gives an identification of $\tilde{\cal G}$ with
${\cal E}nd\bigl(\pi_{_{F*}}\Theta(T^*F)^{^{-1}}\bigl)$
over $E^{(d)}_{_{gen}}$. This, in virtue of the local
freeness of
${\cal E}nd\bigl(\pi_{_{F*}}\Theta(T^*F)^{^{-1}}\bigl)$,
implies our claim.

By definition, $\tilde{\cal G}$ injects in ${\cal R}$. Consider the diagram
$$\matrix{&{\cal R}&\cr
\qquad\nearrow &&\nwarrow\qquad\cr
\widehat{\cal G}&&\tilde{\cal G}\cr
\big\uparrow&&\big\uparrow\cr
\widehat{\cal G}_{\leq A_{ij}({\bf v})} &\buildrel \Phi_{ij{\bf v}}
\over\longleftarrow &\tilde{\cal G}_{\leq A_{ij}({\bf v})} }
\leqno {\bf (4.6.3)}$$
which is commutative by Proposition 4.5 (c). To show the existence
of $\Phi$, it is enough to verify that the $\Phi_{ij{\bf v}}$
preserve the relations between sections of
$\tilde{\cal G}_{\leq A_{ij}({\bf v})}$ which (see Proposition 3.6 (b)),
generate $\tilde {\cal G}$ as an algebra. More precisely, we need to prove
that for any sections $s_k $ of $\tilde{\cal G}_{\leq A_{i_k, j_k}
({\bf v}_k)}$, $k=1, ..., r$, and any relation $F(s_1, ..., s_r)=0$
holding in $\tilde {\cal G}$, we have $F(\Phi(s_1), ..., \Phi(s_r))=0$
(we have omitted the subscripts in the $\Phi$'s). But this is
a consequence of the commutativity of the diagrams (4.6.3)
and of the injectivity of the map $\tilde {\cal G}\rightarrow {\cal R}$.
The compatibility of $\Phi$ with the filtrations is immediate since
the filtrations on $\tilde{\cal G}$ and $\widehat{\cal G}$ are
compatible with the product. Lemma 4.6.2 is proved.

\proclaim (4.6.4) Lemma. For any $A\in {\bf M}$ the homomorphism $\Phi$
induces isomorphisms
$\tilde{\cal G}_{\leq A}/\tilde{\cal G}_{<A}
\buildrel\sim\over\rightarrow
\widehat{\cal G}_{\leq A}/\widehat{\cal G}_{<A}$.

\noindent {\sl Proof.}
As in the proof of Lemma 4.6.2, (2.2.1) and (2.8.8.1) imply that
$\widehat{\cal G}_{\leq A}/\widehat{\cal G}_{<A}\simeq{\cal O}_{E^{(A)}}$.
By construction,
$\tilde{\cal G}_{\leq A}/\tilde{\cal G}_{<A}={\cal O}_{E^{(A)}}$.
Thus, the induced map
$\Phi_A\,:\,\tilde{\cal G}_{\leq A}/\tilde{\cal G}_{<A}
\buildrel\sim\over\rightarrow\widehat{\cal G}_{\leq A}/\widehat{\cal G}_{<A}$
is identified with an endomorphism of the simple sheaf ${\cal O}_{E^{(A)}}$.
Moreover, $\Phi_A$ is non-zero since by localisation of correspondences
it is an isomorphism over $E^{(d)}_{_{reg}}\i E^{(d)}$
(see the proof of (4.6.1)). Thus the map $\Phi_A$ is bijective.

\vskip .3cm

Lemma 4.6.4 implies that $\Phi$ is an isomorphism.

\vskip 2cm

\centerline {\bf 5. Classical elliptic algebras and elliptic cohomology.}

\vskip 1cm

\noindent {\bf (5.1) Heisenberg groups.}
We start with recalling some  well known facts [Mum].
Let $E\to S$ be an elliptic curve over a base scheme $S$.
For any  $n\geq 1$ we denote by $n_{_E}: E\rightarrow E$ the homomorphism of
multiplication by $n$.
 Let $E_n = {\rm Ker} \,n_{_E}\i E$ be the group of points of order $n$, and
let
$\mu_n \i {\bf C}^*$ be the group of $n$th roots of 1.
There is a canonical non-degenerate skew-symmetric pairing
 $$(\alpha,\beta)\mapsto \langle\alpha, \beta\rangle,
\quad E_n \otimes_{\bf Z} E_n
\rightarrow \mu_n$$
called the Weil pairing.
We will sometimes use the notation $\langle\alpha,\beta\rangle_n$
to indicate the dependence on $n$. For instance, if $m=nc$ is a multiple of
$n$, then $E_n\i E_m$, and for any $\alpha,\beta\in E_n$ we have
$\langle\alpha,\beta\rangle_m = \langle\alpha,\beta\rangle_n^c$.

Fix $n\geq 1$. The pairing $\langle\alpha,\beta\rangle$ defines a central
extension
$$1\rightarrow \mu_n \rightarrow H_n \rightarrow E_n\rightarrow 1$$
known as the Heisenberg group. As a set, $H_n = E_n \times \mu_n$, with
multiplication $(\alpha,\zeta) \cdot (\beta, \xi) = (\alpha+\beta,
\, \langle\alpha,\beta\rangle \zeta\xi)$.

\vskip .3cm

\noindent {\bf (5.2) The Heisenberg representation.}
A representation $T$ of $H$ in a vector space (or, more generally,
an action on any variety acted upon by ${\bf C}^*$)
is said to have central charge $c$, if for $\zeta\in\mu_n$ we have
$T(\zeta) = \zeta^c\cdot {\rm Id}$.
For any $c\in {\bf Z}/n$ relatively prime to
$n$ there is a unique, up to isomorphism, $n$-dimensional irreducible
representation  $H_n$ with central charge $c$.
We denote it $T_c: H_n \rightarrow {\rm Aut}(M_c)$.
A particular model for $M_c$ can be obtained as follows. Choose a subgroup
$\Gamma\i E_n$, $\Gamma\simeq {\bf Z}/n$ and define the space
$$M_c(\Gamma) = \bigl\{
f: E_n\rightarrow {\bf C}\,|\, f(\alpha+\gamma) = \langle\alpha,\gamma\rangle^c
f(\alpha)\quad \forall
\alpha\in E_n, \gamma\in\Gamma\bigl\}.$$
 The action of $(\beta,\zeta)\in H_n$ on
$M_c(\Gamma)$ is given by $T_c(\beta, \zeta)(f)(\alpha) =
\zeta^c \langle\alpha,\beta\rangle f(\alpha+\beta)$.

If we choose a complex uniformization $E = {\bf C}/ {\bf Z}\oplus {\bf Z}\tau$,
${\rm Im} (\tau) > 0$, then $E_n$ is identified with
$({\bf Z}/n)^2$ and the
Weil pairing has the form
$$\langle\alpha,\beta\rangle = \exp\biggl( {2\pi i\over n} (\alpha_1\beta_2 -
\alpha_2\beta_1)\biggl), \quad \alpha_i, \beta_i \in {\bf Z}/n.
\leqno {\bf (5.2.1)}$$
The space $M_c(\Gamma)$ for $\Gamma = {\bf Z}/n \oplus 0$, is identified with
${\bf C}^n$ and  the representation $T_c$ sends $(\alpha, \zeta)\in H_n$
into $\zeta^c I_{10}^{c\alpha_1} I_{01}^{\alpha_2}$, where

$$
I_{_{10}}=
\pmatrix {1&0&0&\dots&0\cr0&\omega&0&\dots&0\cr
0&0&\omega^{^2}&\dots&0\cr
\vdots&\vdots&\vdots&\ddots&\vdots\cr
0&0&0&\dots&\omega^{^{n-1}} },\quad
I_{_{01}}=
\pmatrix {0&1&0&\dots&0 \cr 0&0&1&\dots&0 \cr
\vdots&\vdots&\vdots&\ddots&\vdots\cr
0&0&0&\dots&1 \cr 1&0&0&\dots&0}, \quad \omega = e^{2\pi i\over n}.
\leqno{\bf (5.2.2)}$$

\vskip .2cm
For $x\in E$ we denote $t_x: E\rightarrow E$ the translation by $x$.
The main reason for using Heisenberg groups is the following :

\proclaim (5.3) Proposition. Let $L$ be a line bundle of degree
$n$ on $E$. Then for any $\alpha\in E_n$ we have $t_\alpha^*L
\simeq L$. Moreover,
there is an action of $H_n$ in (the total space of)
$L$ extending the action of $E_n$ on $E$ by translations. This action
has central charge 1. The resulting representation of $H_n$ on the space
$\Gamma(E,L)$ is irreducible, i.e. $\Gamma(E,L)\simeq M_1$.

By comparing Weil pairings on $E_n$ and $E_m$, $n|m$, we deduce the
next statement :

\proclaim (5.4) Corollary. If $L$ is a line bundle on $E$ of degree
$nc$, then, for any $E_n$-invariant open set $U\i E$, the space
$\Gamma(U,L)$ is a representation (possibly reducible) of
$H_n$ with central charge $c$.

\vskip .2cm

\noindent {\bf (5.5) Vector bundles over an elliptic curve.}
Let us fix $c$ such that $(c,n)=1$. Let ${\cal M}_E(n,c)$ be the set of
indecomposable rank $n$ vector bundles on $E$ with $c_1(V) = c$.
It was proved by Atiyah [A1] that :

\vskip .2cm

\noindent {\bf (5.5.1)}
 Any two bundles $V,V'\in {\cal M}_E(n,c)$ are related by $V' = V\otimes L$
where $L$ is a line bundle of degree 0.

\vskip .2cm

\noindent {\bf (5.5.2)}
 Any $V\in {\cal M}_E(n,c)$ is simple, i.e. $H^0(E, {\cal E}nd (V)) =
{\bf C}$. Moreover, $H^1(E, {\cal E}nd (V))=0$.

\vskip .2cm

Let us recall, for future use, three equivalent constructions of bundles
from ${\cal M}_E(n,c)$. Each of them gives all the bundles.

\vskip .3cm

\noindent {\bf (5.6) The Heisenberg construction.} Choose a line bundle
$L$ on $E$ of degree $nc$. By Corollary 5.4, the direct image
${n_{_E}}_*L$ has an action of $H_n$ in each fiber, with central charge
$c$. Put $V = {\rm Hom}_{H_n}({n_{_E}}_*L, M_c)$ to be the vector
bundle of spaces of multiplicities of $M_c$ in these representations. Then
$V\in {\cal M}_E(n,c)$.

\vskip .3cm

\noindent {\bf (5.7) The direct image construction.}
Let $\Gamma\i E_n$ be a subgroup isomorphic to ${\bf Z}/n$.
We have a commutative diagram
$$\matrix{ E&\buildrel q\over\longrightarrow &E/\Gamma&\cr
n_{_E}&\searrow &\downarrow&p\cr
&&E&}\leqno {\bf (5.7.1)}$$
where $p$ is an isogeny of degree $n$ with ${\rm Ker} (p) = E_n/\Gamma$.
Then any $V\in {\cal M}_E(n,c)$ can be obtained as $p_*L$ where
$L$ is a line bundle on $E/\Gamma$ of degree $c$.

Note that $\tilde L = q^*L$ is a line bundle on $E$ of degree $nc$,
and $q_*\tilde L$ is canonically identified with $L\otimes M_c(\Gamma)$.
This shows that
$${n_{_E}}_*\tilde L = p_*q_* \tilde L = p_*L\otimes M_c(\Gamma)$$
and thus the constructions (5.7) and (5.6) are equivalent.

\vskip .3cm

\noindent {\bf (5.8) The Fourier transform construction.} Let $P$
 be the Poincar\'e line bundle on $E\times E$ (it corresponds
to the divisor $\Delta - (E\times\{0\}) - (\{0\}\times E)$
where $\Delta$ is the diagonal). Let $p_1, p_2: E\times E \rightarrow E$ be
the projections. Then any $V\in {\cal M}_E(n,c)$ can be obtained as
$$V = R^0 p_{2*} \bigl(p_1^*L\otimes P^{\otimes c}\bigl),
\leqno {\bf (5.8.1)}$$
for a line bundle $L$ on $E$ of degree $n$. For the case $c=1$
this construction is precisely the Fourier-Mukai transform of $L$,
see [Muk].

\vskip .2cm

The construction (5.7) can be found in [Od]. The equivalence of
(5.7) and (5.6) was already explained. To see the equivalence with (5.8),
note that
each fiber of the RHS of (5.8.1) is the space of sections of a line bundle
of degree $n$ on $E$ so has a Heisenberg group action. Taking eigenspaces
of a generator of $\Gamma\i E_n$ we split the fiber into
a direct sum of 1-dimensional vector subspaces. This means that $V$ is the
pushdown of a line bundle from an unramified covering.
We leave the details to the reader.

\vskip .2cm

By (5.5.1), the sheaf
${\cal E}nd(V)$ for $V\in {\cal M}_E(n,c)$ depends only on $n,c$.
We denote this sheaf by ${\cal G}_{n,c}$ or simply $\cal G$.

\proclaim (5.9) Proposition. (a) As a vector bundle, ${\cal G}_{n,c} =
\bigoplus_{\alpha\in E_n} L_\alpha$, where $L_\alpha = {\cal O}_E(0-\alpha)$
is the line bundle of degree 0 corresponding to $\alpha$. The multiplication
$m: {\cal G}_{n,c}\otimes {\cal G}_{n,c}\rightarrow {\cal G}_{n,c}$
restricts to isomorphisms $m_{\alpha \beta}: L_\alpha\otimes L_\beta
\rightarrow L_{\alpha+\beta}$ such that $m_{\beta\alpha} =
\langle\alpha,\beta\rangle^c
\cdot m_{\alpha\beta}$. \hfill\break
(b) The pullback $n_{_E}^*{\cal G}_{n,c}$ (where $n_{_E}:E\rightarrow E$
 is the multiplication by $n$) is identified, as a sheaf of algebras,
with ${\rm Mat}_n({\cal O}_E)$. More precisely, for an open $U\i E$ the
space $\Gamma(U, {\cal G}_{n,c})$ is identified with the space of
matrix functions $A(x) \in {\rm Mat}_n({\cal O}(n_E^{-1}(U)))$
satisfying the condition
$$A(x+\alpha) = T_c(\alpha) A(x) T_c(\alpha)^{-1}, \quad \forall
 \alpha\in E_n.\leqno {\bf (5.9.1)}$$

Part (a) can be found in [A1]. Part (b) follows at once from
the Heisenberg construction of $V$.
 We will call (5.9.1) the automorphy condition, see [RS].

\vskip .3cm

\noindent {\bf (5.10) Classical elliptic algebras.}
Fix a finite subset $P\i E$. The classical elliptic algebra ${\bf e}{\bf
l}_{n,c}$
is, by definition, the Lie algebra of global sections $\Gamma(E\setminus P,
{\cal G}_{n,c})$, where ${\cal G}_{n,c} = {\cal E}nd (E)$ for any bundle
$V\in {\cal M}_E(n,c)$. Let also ${\cal SG}_{n,c} = {\bf s}{\bf l}(V)$ be the
subsheaf
of traceless endomorphisms and ${\bf s}{\bf e}{\bf l}_{n,c}$ be the Lie algebra
of its
global sections over $E\setminus P$. By (5.5.2) the sheaf ${\cal SG}_{n,c}$
has both $H^0$ and $H^1$ vanishing. This makes the subalgebra
${\bf s}{\bf e}{\bf l}_{n,c}$ into a Lie bialgebra by the procedure
described in [Dr1],
[C].
The Lie bialgebra structure is the map
$$\delta: {\bf s}{\bf e}{\bf l}_{n,c} \rightarrow \wedge^2 {\bf s}{\bf e}{\bf
l}_{n,c} \i
{\bf s}{\bf e}{\bf l}_{n,c}\otimes {\bf s}{\bf e}{\bf l}_{n,c} , \quad
\delta(x) = [x\otimes 1 + 1\otimes x, r_{n,c}],$$
where $r_{n,c}\in {\bf s}{\bf e}{\bf l}_{n,c}\hat\otimes {\bf s}{\bf e}{\bf
l}_{n,c}$ is the
classical elliptic $r$-matrix ($\hat\otimes$ means completed tensor product).
Let us describe $r_{n,c}$ explicitly. In the language of Proposition
5.9 (b), this is the meromorphic function $r_{n,c}: E\times E \rightarrow
gl_n \otimes gl_n$ given by
$$r_{n,c}(u,v) = \sum_{\alpha\in E_n - \{0\}} w_\alpha(u-v) T_c(\alpha)
\otimes T_c(\alpha)^{-1},$$
where $w_\alpha$, $\alpha\in E_n$  are meromorphic functions on $E$
uniquely characterized by the following conditions (see [B] for an explicit
formula) :

\item{(1)} $w_\alpha(u+\beta) = {\langle \beta, \alpha\rangle}
w_\alpha(u), \quad \forall \alpha,\beta\in E_n$.

\item{(2)} $w_0\cong 1$ and if $\alpha\neq 0$, then $w_\alpha$ is
regular on $E\setminus E_n$ with simple poles at points of $E_n$ and
${\rm res}_{_{u=0}}w_\alpha(u) = 1$.

\noindent When $c=1$, the matrix $r_{n,1}$ was found by Belavin [B].
The possibility of considering any $c$ prime to $n$ was pointed out
by Cherednik.

\vskip .3cm

\noindent {\bf (5.11) Automorphy   conditions and correspondences.}
We now proceed to construct a realization of the sheaf of  Lie algebras
${\cal G} = {\cal G}_{n,c}$ by correspondences, by using the
description of ${\cal G}$ given by Proposition 5.9 (b).
We suppose that there is a fixed level $n$ structure on $E$,
i.e. an identification
$E_n\buildrel\sim\over\to({\bf Z}/n)^{^2}\times S$ such that the
Weil pairing has the form (5.2.1). Thus the representation
$T_c$ is given by the matrices (5.2.2). We denote $\Gamma={\bf Z}/n\oplus 0$.

We assume the framework and notation of (2.7) (specialized to our case).
Let us identify the set $[n]$ with ${\bf Z}/n$, so ${\bf Z}/n$ acts on $[n]$
by translation. Consider the action of $E_n$ on $E\times [n]^2$ given by
$$(\alpha_1, \alpha_2) \cdot (x, b,c) = (x+\alpha_1, b+\alpha_2, c-\alpha_2).$$
We extend this action  first to
the variety
$$\tilde M= (E\times [n]^2)^{(d)} = \bigoplus_{A\in {\bf M}}
E^{(A)}.$$
Second, we extend the action to an action on the sheaf ${\cal O}_{\tilde
M}$ by setting
$$(\alpha\cdot f)(z) = \omega^{\sum_{i,j}c\alpha_1 a_{ij}(i-j)}\,
f(\alpha\cdot z), \quad z\in \tilde M.\leqno {\bf (5.12)}$$

Let
$${\rm nor}: \tilde M \rightarrow M = (E\times [n])^{(d)}
\times_{E^{(d)}}
(E\times [n])^{(d)}$$ be the natural projection (normalization morphism, see
(2.7)).

Recall (Theorem 3.3) that the variety $M$ is nothing but $Z_{{\cal X}_{U(d)}}$,
the spectrum of the $U(d)$ -equivariant elliptic cohomology of the
Steinberg variety $Z$, and the sheaf ${\rm nor}_* {\cal O}_{\tilde M}$ is
the same as $\Xi_{_Z}$, the microlocal Thom sheaf of the second projection
$Z\rightarrow F$. Recall also
that $\pi_{_Z}: Z_{{\cal X}_{U(d)}}\rightarrow {\cal X}_{U(d)}= E^{(d)}$
is the natural projection.

Set $U = n_{_E}^{-1} (E\setminus P)$. Formula (5.12)
defines an action of $E_n$ in the space $\Gamma(
U^{(d)},{\pi_{_Z}}_*\Xi_{_Z})$.
In the previous section we have constructed a homomorphism of
Lie algebras $\Psi: \Gamma(U, gl_n({\cal O}))\rightarrow
\Gamma( U^{(d)}, {\pi_{_Z}}_*\Xi_{_Z})$.

\proclaim (5.13) Theorem. The map $\Psi$ restricts to a surjective algebra
homomorphism
$${\cal U}({\bf e}{\bf l}_{n,c})
\rightarrow \Gamma( U^{(d)}, {\pi_{_Z}}_*\Xi_{_Z})^{E_n}.$$

\vskip 2cm
\centerline {\bf References}

\vskip 1cm\hangindent=0.7in
[AHS] M. Ando, M.J. Hopkins, N. Strickland, A useful orientation,
Preprint, 1994.

\hangindent=0.7in
[AP]$\enspace$ C. Allday, V. Puppe, Cohomological Methods in Transformation
Groups (Cambridge Studies in Math., {\bf 32}), Cambridge University
Press, 1993.

\hangindent=0.7in
[A1]$\enspace$ M. Atiyah, Vector bundles over an elliptic curve, {\it Proc.
London
Math. Soc.} {\bf 7} (1957), 414-452.

\hangindent=0.7in
[A2] $\enspace$M. Atiyah, K-theory, Benjamin, 1968.

\hangindent=0.7in [BLM] A.A. Beilinson, G. Lusztig, R.W. McPherson,
A geometric setting for the quantum deformation of $GL_n$,
{\it Duke Math. J.}, {\bf 61} (1990), 655-675.

\hangindent=0.7in
[B]$\enspace\enspace$ A.A. Belavin, Discrete groups and integrability of
quantum systems,
{\it Funct. Anal. Appl.} {\bf 14} (1981), 260-267.

\hangindent=0.7in
[BD]$\enspace$ A.A. Belavin, V.G. Drinfeld, Solutions of classical Yang-Baxter
equations
for simple Lie algebras, {\it Funct. Anal. Appl.} {\bf 16} (1983),
159-180.

\hangindent=0.7in
[BG] $\enspace$J. Block, E. Getzler, Equivariant cyclic homology and
equivariant
differential forms.

\hangindent=0.7in [BS] $\enspace$I.N. Bernstein, O.V. Schwartzman, Chevalley
theorem for
complex crystallographic groups,
{\it Funct. Anal. Appl.}, {\bf 12}:4 (1978), 79-80.

\hangindent=0.7in [C]$\enspace\enspace$ I. V. Cherednik, Algebraic Methods in
Soliton Theory, Moscow, 1994
(in Russian).

\hangindent=0.7in
[CG] $\enspace$N. Chriss, V. Ginzburg, Representation theory and Complex
Geometry
(Geometric technique in Representation theory of reductive groups),
Progress in Mathem., Birkhauser 1995.

\hangindent=0.7in
[Di]$\enspace$ T. tom Dieck, Lokalizierung \"aquivarianter
Kohomologie-Theorien,
{\it Math. Z.} {\bf 121} (1971), 253-262.

\hangindent=0.7in [Dr1]$\enspace$ V. Drinfeld, Quantum groups, Proc. ICM-86,
Berkeley.

\hangindent=0.7in[Dr2]$\enspace$ V. Drinfeld, Yangians and degenerate affine
Hecke algebras,

\hangindent=0.7in
[EFK] P. Etingof, I. Frenkel, A. Kirillov, Jr. Spherical functions
on affine Lie groups, preprint, 1994 (hep-th 9407047).

\hangindent=0.7in [GKV] V. Ginzburg, M. Kapranov, and E. Vasserot.
Langlands reciprocity for algebraic surfaces. {\it Mathem. Research
Letters},
{\bf 2} (1995), 147-160.

\hangindent=0.7in [Gr]$\enspace$  I. Grojnowski. Delocalized equivariant
elliptic cohomology.
      Yale University preprint, 1994.

\hangindent=0.7in
[GV1]  V. Ginzburg et E. Vasserot, Langlands Reciprocity for Affine Quantum
groups of type $A_n$, {\it Intern. Mathem. Research Notices. Duke
Math. J.}, {\bf 71}, 1993, p. 67-85.

\hangindent=0.7in
[GV2] V. Ginzburg et E. Vasserot. Alg\`ebres elliptiques et $K$-th\'eorie
\'equivariante, {\it C.R. Acad.Sci. Paris}, {\bf 319} (1994), 539.

\hangindent=0.7in
[HBJ] F. Hirzebruch, Th. Berger, R. Jung, Manifolds and Modular Forms
(Aspects of Mathematics, vol. E-20), F. Vieweg Publ., Braunschweig 1992.

\hangindent=0.7in
[KL]$\enspace$ D. Kazhdan, G. Lusztig, Proof of Deligne-Langlands conjecture
for
affine Hecke algebras, {\it Invent. Math.}, {\bf 87} (1987), 153-215.

\hangindent=0.7in
[Lan] P. Landweber (Ed.) Elliptic Curves and Modular Forms in Algebraic
Topology (Lecture Notes in Math., {\bf 1326}), Springer-Verlag, 1986.

\hangindent=0.7in
[Muk] S. Mukai, A duality between $D(E)$ and $D(E^\vee)$,
 {\it Nagoya Math. Journal}

\hangindent=0.7in
[Mum] D. Mumford, Tata Lectures on Theta I (Progress in Math. {\bf 28}),
Birkhauser, 1983.

\hangindent=0.7in
[Q]$\enspace\enspace$ D. Quillen, Elementary proof of some properties of
complex cobordism
via Steenrod operations, {\it Adv. in Math.} {\bf 7} (1972), 29-56.

\hangindent=0.7in
[Od]$\enspace$ T. Oda, Vector bundles on an elliptic curve, {\it Nagoya Math.
J.}
{\bf 43} (1971), 41-72.

\hangindent=0.7in
[R]$\enspace\enspace$ A. Ramanathan, Stable Principal Bundles on a Compact
Riemann
Surface, {\it Math. Ann.}, {\bf 213} (1975), 129-152.

\hangindent=0.7in
[RS]$\enspace$ A.G. Reiman, M.A. Semenov -Tyan-Shanskii, {\it Journ. Soviet.
Math.},
{\bf 46} (1989), 1631.

\hangindent=0.7in
[S1]$\enspace$ G. Segal, Equivariant K-theory, {\it Publ. IHES}, {\bf 34}
(1968), 129-151.

\hangindent=0.7in
[S2]$\enspace$ G. Segal, Elliptic cohomology, {\it Expos\'e au S\'eminaire
Bourbaki.}
{\bf 695} (1988).

\hangindent=0.7in
[V]$\enspace\enspace$ E. Vasserot, Repr\'esentations de groupes quantiques and
permutations,
{\it Ann. ENS}, {\bf 26} (1993), 747-773.
\vskip 2cm

V.G.: The University of Chicago, Mathem. Dept., Chicago, IL 60637, USA
\vskip 5mm

M.K.: Northwestern University, Mathem. Dept., Evanston, IL 60208,
USA
\vskip 5mm

E.V.: Ecole Normale Sup\'erieure, 45 rue d'Ulm, 75005 Paris, FRANCE

\bye